\documentclass{article}

 \usepackage{amssymb,amsmath}

\begin{document}

\def\Q{{\Bbb Q}}
 \def\R{{\mathbb R}}
 \def\C{{\mathbb C}}
\def\O{{\Bbb O}}
\def\K{{\Bbb K}}
\def\L{{\rm Lat}}
\def\Z{{\mathbb Z}}
\def\F{{\Bbb F}}
\def\Q{{\Bbb Q}}

\def\cV{{\cal V}}
\def\cH{{\cal H}}
\def\cA{{\cal A}}

\def\fP{{\frak P}}
\def\fI{{\frak I}}
\def\fN{{\frak N}}
\def\fF{{\frak F}}

\def\kk{{\bold k}}

\def\phi{\varphi}
\def\kappa{\varkappa}
\def\epsilon{\varepsilon}
\def\le{\leqslant}
\def\ge{\geqslant}

\def\GL{{\rm GL}}
\def\Sp{{\rm Sp}}

\def\Fl{{\rm Fl}}

\def\wt{\widetilde}

\newcommand{\codim}{\mathop{\rm codim}\nolimits}
\newcommand{\vol}{\mathop{\rm vol}\nolimits}
\newcommand{\res}{\mathop{\rm res}\nolimits}
\renewcommand{\Re}{\mathop{\rm Re}\nolimits}

\renewcommand{\Im}{\mathop{\rm Im}\nolimits}

\newcounter{sec}
\renewcommand{\theequation}{\arabic{sec}.\arabic{equation}}

\def\2{
\begin{picture}  (300,130)(0,-20)
\def\bm{\begin{matrix}}
\def\em{\end{matrix}}
\def\ss{\scriptstyle}
\put(0,0){\vector(1,0){320}}

\multiput(10,0)(0,5){17}{\line(0,1){3}}
\put(12,7){$I_3,\,I_{21}^\pm,\,I_{12}^\pm,\,I_{03}^\pm$}
\put(12,25){$I_2,\,I_{11}^\pm,\,I_{02}^\pm$}
\put(12,45){$I_1,\,I_{01}^\pm$}
\put(12,65){$I_0$}
\put(5,-15){$I_{31}^\pm$}
\put(45,-15){$I_{3}$}

\put(0,90){$\bm\ss {\rm Re}\,\alpha=\\ \ss =n-7\em$}

\multiput(90,0)(0,5){17}{\line(0,1){3}}
\put(92,25){$I_2,\,I_{11}^\pm,\,I_{02}^\pm$}
\put(92,45){$I_1,\,I_{01}^\pm$}
\put(92,65){$I_0$}
\put(85,-15){$I_{21}^\pm$}
\put(125,-15){$I_{2}$}

\put(80,90){$\bm\ss {\rm Re}\,\alpha=\\ \ss =n-5\em$}

\multiput(170,0)(0,5){17}{\line(0,1){3}}
\put(172,45){$I_1,\,I_{01}^\pm$}
\put(172,65){$I_0$}
\put(165,-15){$I_{11}^\pm$}
\put(205,-15){$I_{1}$}

\put(160,90){$\bm\ss {\rm Re}\,\alpha=\\ \ss =n-3\em$}

\multiput(250,0)(0,5){17}{\line(0,1){3}}
\put(252,65){$I_0$}
\put(245,-15){$I_{01}^\pm$}
\put(240,90){$\bm\ss {\rm Re}\,\alpha=\\ \ss =n-1\em$}


\multiput(10,0)(40,0){7}{\circle*{3}}

\put(290,-4){\line(0,1){8}}
\put(287,6){$\scriptstyle n$}

\put(315,6){$\alpha$}

\end{picture}

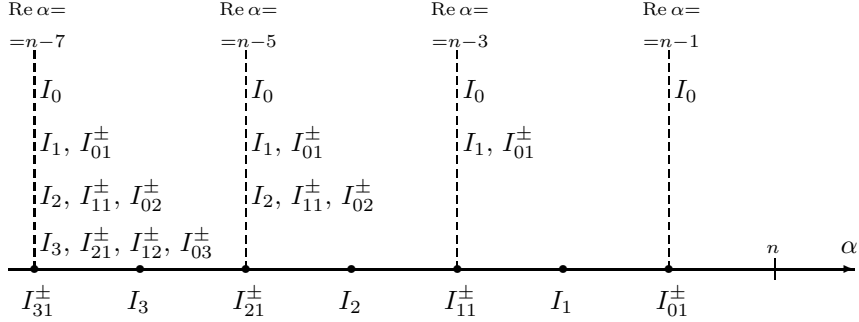
{\sf Fig 2. The summands of the Plancherel formula in
different  strips. The summand of the Plancherel formula surviving at
the points $\alpha=n-1-m\ge 0$.}
       }

\def\1{
{\sf Fig 1.}

\begin{picture} (300,80) (-10,-40)
\put(0,-4){\line(0,1){8}}
\put(-2,5){0}
\put(0,0){\vector(1,0){300}}
\qbezier(80,20)(85,0)(80,-20)
\qbezier[60](160,40)(170,0)(160,-40)
\qbezier[20](40,10)(42.5,0)(40,-10)

\put(80,25){$\scriptstyle |z_1|=1$}
\put(170,25){$\scriptstyle |z_1|=p$}
\put(30,25){$\scriptstyle |z_1|=1/p$}

{
\thicklines
\put(95,0){\vector(-1,0){8}}
\put(95,0){\circle*{3}}
\put(73,0){\vector(1,0){8}}
\put(73,0){\circle*{3}}
\put(95,5){$\scriptstyle p^{\beta}$}
\put(67,5){$\scriptstyle p^{-\beta}$}
}

\end{picture}

{\sf a) The poles $z_1=p^{\pm\beta}$ of the integrand $I_0$
for $\beta>0$
and their motion for decreasing $\beta$.
If $z_2$,\dots, $z_n$ are fixed and $|z_2|=\dots=|z_n|=1$,
then all other poles (with respect $z_1$) lie on the circles
$|z_1|=p^{-1}$, $|z_1|=p$.}

\begin{picture} (300,80) (-10,-40)
\put(0,-4){\line(0,1){8}}
\put(0,0){\vector(1,0){320}}
\qbezier(80,20)(85,0)(80,-20)
\qbezier[60](160,40)(170,0)(160,-40)
\qbezier[20](40,10)(42.5,0)(40,-10)


{
\thicklines
\put(95,0){\vector(-1,0){8}}
\put(95,0){\circle*{3}}
\put(73,0){\vector(1,0){8}}
\put(73,0){\circle*{3}}
\put(95,5){$\scriptstyle p^{\beta+1}$}
\put(60,5){$\scriptstyle p^{-\beta-1}$}
\put(290,0) {\vector(1,0){8}}
\put(290,0) {\circle*{3}}
\put(285,5){$\scriptstyle p^{-\beta+1}$}
\put(10,0){\vector(-1,0){8}}
\put(10,0){\circle*{3}}
\put(5,5){$\scriptstyle p^{\beta-1}$}
}

\end{picture}

{\sf b) The poles of the integrand $I_1$ for
$0>\beta>-1$. The $\beta$ decreases.  }

\begin{picture} (300,80) (-10,-40)
\put(0,-4){\line(0,1){8}}
\put(0,0){\vector(1,0){320}}
\qbezier(80,20)(85,0)(80,-20)
\qbezier[60](160,40)(170,0)(160,-40)
\qbezier[20](40,10)(42.5,0)(40,-10)


{
\thicklines
\put(73,0){\vector(1,0){8}}
\put(73,0){\circle*{3}}
\put(57,5){$\scriptstyle p^{-\beta-1}$}
\put(290,0) {\vector(1,0){8}}
\put(290,0) {\circle*{3}}
\put(285,5){$\scriptstyle p^{-\beta+1}$}
\put(45,0){\vector(-1,0){8}}
\put(45,0){\circle*{3}}
\put(42,5){$\scriptstyle p^{\beta}$}
}

\end{picture}

{\sf c) The poles of the integrand $I_{01}^-$ for
$0>\beta>-1$. The $\beta$ decreases. }
                                  }

\begin{center}

{\large\bf

Beta-function of  Bruhat--Tits buildings
and  deformation of $l^2$ on the set of
$p$-adic lattices
}

\medskip

\sc\large Yurii A. Neretin\footnote
{partially supported by the
grant
NWO 047-008-009}

\end{center}

{\small
For  the space $\L_n$ of all the lattices
in a $p$-adic $n$-dimensional linear space we obtain an analog
of matrix beta-function; this beta-function
has a degeneration to the Tamagawa zeta-function.
We propose an analog of Berezin kernels
for $\L_n$. We obtain conditions of positive definiteness
of these kernels  and explicit Plancherel formula.
       }

\medskip

It is well known that affine Bruhat--Tits buildings
(see, for instance, \cite{Bro}, \cite{Gar})
are right $p$-adic analogs of Riemannian noncompact
symmetric spaces.  This analogy exists
on the levels of  geometry,  harmonic analysis,
and special functions.
We continue this parallel.

First, we obtain an imitation of a matrix beta-function
(see \cite{Gin}, \cite{Ner1}) for the spaces
\begin{equation}
\GL_n(\Q_p)/\GL_n(\Z_p)
,\end{equation}
see Theorem 2.1, this beta extends the Tamagawa zeta
(Corollary 2.12, Theorem 2.13).

Second, we define an analog of Berezin kernels
(see \cite{Ner2})
for the homogeneous spaces (0.1).
We obtain  conditions of positive definiteness of these kernels
(Theorem 2.2)
and the Plancherel formula (Theorems 2.3--2.4).

It seems, that there is no natural $p$-adic analogs
of Hermitian Cartan domains
and holomorphic discrete series.
Nevertheless, the Berezin kernels on Riemannian symmetric spaces
  admit an  imitation for buildings.
More precisely, our objects are spaces of lattices,
i.e., sets of vertices of buildings.
 These 'Berezin kernels'    (see formula (2.5))
have sense for any  building associated with
a classical $p$-adic group;
in this paper,
we consider only the series of groups $\GL_n$.
In this case the Plancherel formula, the conditions of positive
definiteness,
limit behavior are similar to the real case
as precisely as it is possible.
By an analogy with the real case,
it seems that for sympletic and especially orthogonal groups
the theory must be more rich.
It is natural to wait for appearence of more complicated
spectra (as in \cite{Ner2}) and of operators of restriction
to some ideal boundary (as in \cite{Ners}) separating these
spectra;
the both phenomena really appear for  Bruhat--Tits trees,
see \cite{Ner3}.

For the Bruhat--Tits trees our kernels
were known earlier (see \cite{Haa}, \cite{Ols}),
in this  case (it is not discussed in this paper)
the   parallel between the Lobachevsky plane
and the tree also can be extended
to the group of diffeomorphisms
of circle, see
 \cite{Ner3}.

Section 1 contains preliminaries,
the results of this work are formulated in
Section 2, proofs are contained in Sections 3--6.

\medskip

{\bf\large  1. Preliminaries}

\addtocounter{sec}{1}
\setcounter{equation}{0}

\nopagebreak

\medskip

\nopagebreak

{\bf 1.1. Notation.}
Denote by $\K$ the field $\Q_p$  of $p$-adic numbers.
By $\O$ we denote the ring of $p$-adic integers.
By $|\cdot|$ we denote
the norm on $\K$. If $z,z^{-1}\in \O$, then
$|z|=1$. Also $|p|=p^{-1}$, and $|zu|=|z|\,|u|$.

By $\F_p$ we denote the field with $p$ elements; we have
$\F_p\simeq\O/p\O$.

Denote by $\K^n$ the $n$-dimensional linear space
over $\K$. Denote by $e_1$, \dots, $e_n$ the standard basis
in $\K^n$.  By $\K^l$, we denote the subspace in $\K^n$
spanned by $e_1$, \dots, $e_l$.

By $\GL_n(\K)$  we denote the group of invertible
$n\times n$ matrices
over $\K$. By $\GL_n(\O)$ we denote the group of matrices
$g\in\GL_n(\K)$ such that all the matrix elements of $g$
and $g^{-1}$ are integer.
The group $\GL_n(\O)$ is a maximal compact subgroup
in $\GL_n(\K)$.

By  $\L_n$ we denote the space of all  lattices
(see a definition in 1.2)
in  $\K^n$.
By $l^2(\L_n)$ we denote the Hilbert space of
complex-valued functions $f$ on $\L_n$
satisfying the condition $\sum_R |f(R)|^2<\infty$;
this space is
equipped with the scalar product
$$
\langle f_1,f_2\rangle=
\sum_{R\in \L_n}
f_1(R)\overline{f_2(R)}
.$$

\smallskip

{\bf 1.2. Spaces of lattices.}
A {\it lattice} in $\K^n$ is a compact open $\O$-submodule.
 Recall, that any lattice $R$ can be represented
in the form
\begin{equation}
R=\O f_1\oplus \dots\oplus\O f_n
,
\end{equation}
where $f_1$, \dots $f_n$ is some basis in $\K^n$, see
\cite{Wei}, \S II.2.
By $\O^n$ we denote the lattice
$\O e_1\oplus \dots\oplus\O e_n$.

The group $\GL_n(\K)$ acts transitively on $\L_n$.
The stabilizer of the lattice
$\O^n$
is $\GL_n(\O)$. Thus,
$$\L_n=\GL_n(\K)/\GL_n(\O).$$

Orbits of the compact group $\GL_n(\O)$
on $\L_n$ are finite. Each orbit contains
a unique representative of the form
\begin{equation}
 p^{k_1}\O e_1\oplus\dots\oplus p^{k_n}\O e_n,
\qquad \text{where $k_1\ge \dots\ge k_n$}
,\end{equation}
see \cite{Wei}, \S II.2, Theorem 2.

In other words, the space $\GL_n(\O)\setminus\L_n$
  of $\GL_n(\O)$-orbits    on $\L_n$
is in one-to-one correspondence with the set $\fP^+_n$,
whose elements are collections of integers
\begin{equation}
\kk:\qquad k_1\ge\dots\ge k_n
.\end{equation}
Denote by  ${\cal O}[\kk]$ the orbit containing
the lattice (1.2).

The set  $\fP^+_n=\GL_n(\O)\setminus\L_n$   coincides
with the space of double cosets
$\GL_n(\O)\setminus\GL_n(\K)/ \GL_n(\O)$.

For  $s\in \Z$, denote by $m_s$ the number of $k_j$ equal to $s$.
The number of points in the orbit  ${\cal O}[\kk]$
is
\begin{equation}
\nu(\kk)=
\frac{p^{\sum (n-2j+1)k_j}\prod_{1\le l\le n}(1-p^{-l}) }
{\prod_s \prod_{1\le l\le  m_s}(1-p^{-l})}
,\end{equation}
see \cite{Mac2}, formula (V.2.9).

\smallskip

{\bf 1.3. Volume.}
There exists a unique up to a scalar factor
translation invariant measure ({\it volume}) on $\K^n$.
We normalize the volume in $\K^n$ by the condition
$$\vol(\O^n)=1.$$
Obviously,
$$\vol(\oplus_{j=1}^n p^{k_j}\O e_j)=p^{-\sum k_j}.$$

We also normalize the volume in any linear subspace
$V\subset \K^n$ by the condition
$$\vol_V(\O^n\cap V)=1.$$

 {\bf 1.4. Spherical functions on $\GL_n(\K)$.}
Denote by $B_n$ the group of all upper triangular
$n\times n$ matrices over $\K$, the group $B_n$ is the stabilizer
of the flag
$$
0=\K^0\subset \K^1 \subset\dots\subset \K^{n-1}\subset \K^n
.$$

Denote by $\Fl_n$ the space of all  flags
$$\cV:\,\,
0=V_0\subset V_1\subset V_2\subset
   \dots \subset V_{n-1}\subset V_n=\K^n;
\qquad \text{where $\dim V_j=j$}
,$$
in $\K^n$; the space $\Fl_n$ is the homogeneous space
$\GL_n(\K)/B_n$. This space
is also  homogeneous with respect
to the compact group $\GL_n(\O)$;
 hence there exists a unique up to a scalar factor
$\GL_n(\O)$-invariant measure $d\cV$  on $\Fl_n$.
We assume that the total measure of $\Fl_n$ is 1.

For an element $V_j$ of a flag $\cV$, consider the natural
measure $\vol_{V_j}$ on $V_j$
and the natural measure $\vol_{gV_j}$  on $gV_j$.
 Consider also the image $g\cdot\vol_{V_j}$ of $\vol_{V_j}$
under the map $g$. Thus we have two measures on $gV_j$,
denote by $a_j(g;\cV)$ their ratio, i.e.,
$$
g\cdot\vol_{V_j}(Q) =a_j(g;\cV) \vol_{gV_j} (Q)
$$
for each subset $Q\subset gV_j$.

Obviously, these numbers (multipliers)
 satisfy the identity
\begin{equation}
a_j(g_1g_2;\cV)= a_j(g_1,g_2\cV)
                         a_j(g_2;\cV)
.\end{equation}
The Radon--Nikodym derivative
$d(g\cV)/d\cV$
of the transformation $\cV\mapsto g\cV$
is
\begin{equation}
\frac{d(g\cV)}{d\cV}=
\prod\limits_{j=1}^{n}
 \Bigl(\frac{a_j(g;\cV)}  {a_{j-1}(g;\cV)}
     \Bigr)^{-(n+1)+2j}
.\end{equation}

Fix complex numbers $\lambda_1,\dots,\lambda_n$.
Consider the representation $T_\lambda$ of the group $\GL_n(\K)$
in the space $L^2(\Fl_n)$  defined by the formula
\begin{multline}
T_\lambda(g) f(\cV)=
f(g \cV)\cdot
\prod\limits_{j=1}^{n-1}
a_j(g;\cV)^{-1-\lambda_j+\lambda_{j+1}}
\cdot a_n(g;\cV)^{(n-1)/2-\lambda_n}
=\\=
f(g\cV)\cdot
\prod\limits_{j=1}^n
 \Bigl(\frac{a_j(g;\cV)}  {a_{j-1}(g;\cV)}
     \Bigr)^{-(n+1)/2+j-\lambda_j}
\end{multline}
(we assume $a_0=1$).  It is a representation, since (1.5).
If $\lambda_j=i s_j$ are pure imaginary, then $T_\lambda$
is unitary in $L^2$ (this follows from (1.6)).

The representations $T_\lambda$  are called  representations
of  the {\it nondegenerate
 principal series}. For imaginary $\lambda_j=is_j$,
they are called  representations
of the {\it nondegenerate  unitary principal series.}

Let $\sigma\in S_n$ be a permutation of $(\lambda_1,\dots,\lambda_n)$.
For pure imaginary $\lambda_j=i s_j$,
we have
$$T_{\sigma\lambda}\simeq T_\lambda,$$
the same is valid for generic $\lambda\in\C^n$.

Denote by $\bold 1$ the function $f(\cV)=1$ on $\Fl_n$.
The {\it spherical function} $\phi_\lambda(g)$ is defined
as the following matrix element of the representation $T_\lambda$
\begin{equation}
\phi_\lambda(g):=
\langle T_\lambda(g){\bold 1},{\bold 1} \rangle_{L^2}
=
\int\limits_{\Fl_n}\prod\limits_{j=1}^n
 \Bigl(\frac{a_j(g,\cV)}  {a_{j-1}(g,\cV)}
     \Bigr)^{-(n+1)/2+j-\lambda_j}   d\cV
,\end{equation}
where $g\in\GL_n(\K)$.

Spherical  functions are $\GL_n(\O)$-biinvariant
in the following sense
\begin{equation}
\phi_\lambda(g)=\phi_\lambda(h_1 g h_2),\qquad
\text{where $g\in\GL_n(\K)$, $h_1$, $h_2\in\GL_n(\O)$}
.\end{equation}
Indeed,
$$
\phi_\lambda(h_1 g h_2)
=\langle T_\lambda(h_1 g h_2){\bold 1},{\bold 1} \rangle
=\langle T_\lambda( g) T_\lambda( h_2){\bold 1},
T_\lambda^*(h_1){\bold 1} \rangle
=\langle T_\lambda(g){\bold 1},{\bold 1} \rangle
.$$

These functions are symmetric with
respect to  permutations of the parameters $\lambda_j$
$$
\phi_{\sigma\lambda}(g)=\phi_\lambda(g),\qquad \sigma\in S_n
.$$

Also for each  vector
\begin{equation}
\upsilon=\tfrac {2\pi i}{\ln p}
\bigl( l_1, \dots, l_n\bigr),
\qquad
l_j\in\Z
,\end{equation}
we have
$$
\phi_{\lambda+\upsilon}(g)=\phi_{\lambda}(g)
.$$
Indeed, the factors $a_j$ in (1.8)
satisfy
$\log_p a_j\in\Z$.

Thus, the spherical   functions $\phi_\lambda$
are parametrized   by the quotient  set
$$
\C^n\bigl/\bigl[S_n\ltimes \tfrac {2\pi i}{\ln p}\,\Z^n\bigr]
$$  of $\C^n$ by
the group generated by permutations of
the coordinates and shifts
by  vectors (1.10).

We have defined a spherical function $\phi_\lambda$
as a  function on $\GL_n(\K)$. By the biinvariance (1.9),
we can consider it as a function on
 $\L_n=\GL_n(\K)/\GL_n(\O)$ or as a function on
the double cosets $\fP_n^+=\GL_n(\O)\setminus\GL_n(\K)/\GL_n(\O)$.

\smallskip

{\bf 1.5. Spherical functions as functions on
$\GL_n(\O)\setminus\GL_n(\K)/\GL_n(\O)$.}
An explicit expression for the spherical functions
as functions on $\fP^+_n$
 was obtained in \cite{Mac1}, see also \cite{Mac2},
formula V.(3.4):
\begin{equation}
\phi_\lambda(k_1,\dots,k_n)=A\cdot
p^{\sum  ((n+1)/2-j) k_j }
\sum\limits_{\sigma \in S_n}
p^{-\sum_j k_j \lambda_{\sigma(j)}}
\prod\limits_{m<l}
\frac{p^{-\lambda_{\sigma(m)}}-p^{-1-\lambda_{\sigma(l)}}}
    {p^{-\lambda_{\sigma(m)}}-p^{-\lambda_{\sigma(l)}}}
,\end{equation}
where the summation is taken over all the permutations
$\sigma$ of the set $\{1,2,\dots,n\}$ and
the constant $A$ is
$$A=\frac{(1-p^{-1})^n}
 {\prod_{j=1}^n (1-p^{-j})}
.$$

{\sc Remark.} Thus, the expression $\phi_\lambda(\kk)$
for a fixed $\kk$ is a Hall--Littlewood  symmetric function
(see \cite{Mac2}, Chapter III)
as a function in the variables $z_j=p^{\lambda_j}$.

\smallskip

{\bf 1.6. Spherical functions as functions on $\L_n$.}
A spherical function
 $\phi_\lambda$ on $\L_n$
has the following integral representation
\begin{equation}
\phi_\lambda(R)=
\int\limits_{\GL_n(\O)} \prod\limits_{j=1}^{n-1}
\vol\bigl((hR)\cap \K^j\bigr)^{1+\lambda_j-\lambda_{j+1}}
\cdot\vol(hR)^{-(n-1)/2+\lambda_n}
dh
,\end{equation}
where $dh$ is the Haar measure
on $\GL_n(\O)$; we assume   that the total measure
of $\GL_n(\O)$ is 1.  This  follows
from (1.9): since the  $\Fl_n$ is a $\GL_n(\O)$-homogeneous
space, we can replace the integration
over $\Fl_n$ by the integration over $\GL_n(\O)$.


Next,  we  replace the integration over $\GL_n(\O)$
in (1.12) by the average over   $\GL_n(\O)$-orbits
\begin{equation}
\phi_\lambda(R)=
\frac1{\nu(\kk)}
\sum\limits_{R\in{\cal O}[\kk]}
\prod\limits_{j=1}^{n-1}
\vol\bigl(R\cap \K^j\bigr)^{1+\lambda_j-\lambda_{j+1}}
\cdot\vol(R)^{-(n-1)/2+\lambda_n}
,\end{equation}
where $\nu(\kk)$ is given by (1.4).

\smallskip

{\bf 1.7. Spherical functions as invariant kernels on $\L_n$.}
Let $L(R,S)$ be a function on $\L_n\times\L_n$ invariant
with respect to the action of $\GL_n(\K)$,
i.e.,
$$L(gR,gS)=L(R,S), \qquad
\text{for $R$, $S\in\L_n$, $g\in\GL_n(\K)$.} $$
Then
$$\ell(S):=L(\O^n,S)$$
is a $\GL_n(\O)$-invariant function on $\L_n$. Conversely, having
a $\GL_n(\O)$-invariant function
$\ell$ on $\L_n$, we can reconstruct the
corresponding $\GL_n(\K)$-invariant kernel
$L(R,S)$.

Let us apply this remark to the spherical functions
and define the {\it spherical kernel} $\Phi_\lambda(R,T)$
on $\L_n\times \L_n$ by the formula
$$
\Phi_\lambda( g \O^n, h\O^n):=\phi_\lambda (g^{-1}h\O^n)
.$$

{\bf 1.8. Spherical transform.}
Let $f$ be a $\GL_n(\O)$-invariant function
on $\L_n$. Its {\it spherical transform} is the function
$\wt f$ in  variables $\lambda_j$ defined by
\begin{equation}
\widehat f(\lambda_1,\dots,  \lambda_n)
=\sum\limits_{R\in\L_n}
\phi_{-\lambda_1,\dots,  -\lambda_n}(R)
f(R)
.\end{equation}

If we consider $f$ as a function on $\fP^+_n$, we write
$$
\widehat f(\lambda_1,\dots,  \lambda_n)
=\sum\limits_{\kk\in\fP^+_n} \nu(\kk)
\phi_{-\lambda_1,\dots,  -\lambda_n}(\kk)
f(\kk)
.$$

Since $f(R)$ in (1.14) is a constant on any $\GL_n(\O)$-orbit,
formula (1.13) allows to convert (1.14) to the form
\begin{equation}
\widehat f(\lambda_1,\dots,  \lambda_n)
=\sum\limits_{R\in\L_n} f(R)    \Bigl[
\prod\limits_{j=1}^{n-1}
\vol\bigl(R\cap \K^j\bigr)^{1+\lambda_j-\lambda_{j+1}}
\cdot\vol(R)^{-(n-1)/2+\lambda_n}
\Bigr]
.\end{equation}

{\bf 1.9. Plancherel theorem.}
 Macdonald's {\it inversion formula} for the spherical transform
(see \cite{Mac1}) is
\begin{equation}
f(R)=\int\limits_{0\le s_j\le 2\pi/\ln p}
\widehat f(is_1,\dots,is_n)
 \phi_{is_1,\dots,is_n}(R)\, d\mu(s)
,\end{equation}
where the {\it  Plancherel measure} $d\mu(s)$  is
\begin{equation}
d\mu(s)=
C\cdot \prod\limits_{k<l}
\Bigl|
\frac{p^{is_k}-p^{is_l}} {p^{is_k}-p^{-1+is_l}}
\Bigr|^2    ds_1\dots ds_n
\end{equation}
and the constant $C$ is
\begin{equation}
C=\frac  {\ln^n p}{ (2\pi)^n n! }
 \prod\limits_{j=1}^n \frac{1-p^{-j}}{1-p^{-1}}
.\end{equation}

Also the following  {\it Plancherel formula}
holds
$$
\sum\limits_{\kk\in\fP^+_n} |f(\kk)|^2\nu(\kk)=
\int\limits_{0\le s_j\le 2\pi/\ln p}
|\widehat f(is)|^2\,d\mu(s)
.$$

{\bf 1.10. Explicit spectral decomposition of $l^2(\L_n)$.}
Denote by $\Xi_n$
the simplex
\begin{equation}
0\le s_1\le\dots\le s_n\le 2\pi/\ln p
.\end{equation}
Equip this simplex by the Plancherel measure
$n!\,d\mu(s)$, where $d\mu(s)$ is given by (1.17).
Consider the space $\Xi_n\times\Fl_n$
equipped with the product-measure $d\mu(s)\times d\cV$.

The group $\GL_n(\K)$ acts
in the space $L^2(\Xi_n\times\Fl_n)$
 by the  unitary operators
\begin{equation}
\rho(g)\Psi(s,\cV)=\Psi(s,g \cV)\cdot
\prod\limits_{j=1}^n
 \Bigl(\frac{a_j(g,\cV)}  {a_{j-1}(g,\cV)}
     \Bigr)^{-(n+1)/2+j-is_j}
.\end{equation}
This formula defines
 the direct multiplicity-free integral of
unitary representations
of principal nondegenerate series.

For a function $f(R)$ on $\L_n$, we define the function
$J f(s,\cV) $ given by
\begin{equation}
Jf(s,\cV)=
\sum\limits_{R\in\L_n}
f(R)\cdot \prod\limits_{j=1}^{n-1}
\vol\bigl(R\cap \K^j\bigr)^{1-is_j+is_{j+1}}
\cdot\vol(R)^{-(n-1)/2-is_n}
.\end{equation}

It can easily be checked, that the operator $J$ intertwines
representation of $\GL_n(\K)$ in $l^2(\L_n)$
and the representation $\rho$.
Also, the operator $J$ is a unitary operator
$$
l^2(\L_n)\to
L^2(\Xi_n\times\Fl_n,d\mu(s)\times d\cV)
$$
(this is a rephrasing of the Plancherel formula for   the
spherical transform).

\medskip

{\bf\large 2. Results}

\addtocounter{sec}{1}
\setcounter{equation}{0}

\nopagebreak

\medskip

{\bf 2.1. Beta-function.}
Consider the  flag
\begin{equation}
0=\K^0\subset \K^1\subset \dots \subset \K^{n-1}\subset \K^n
.
\end{equation}
Consider the lattices
$
\O^j=\K^j\cap \O^n$.
Thus we obtain the flag of $\O$-modules
$$
0=\O^0\subset \O^1\subset \dots \subset \O^{n-1}\subset \O^n
.
$$

Fix complex numbers
$\alpha_1,\dots, \alpha_n, \beta_1, \dots, \beta_n$
such that
\begin{equation}
\Re\beta_j+j-1<0,
\qquad
\Re\alpha_j+\Re\beta_j-n+j>0
.\end{equation}
Let also $\alpha_{n+1}=\beta_{n+1}$.

\smallskip

{\sc Theorem 2.1.}
\begin{gather}
\sum_{R\in \L_n}
\prod\limits_{j=1}^n
\Bigl\{
\vol(R\cap \K^j)^{\beta_j-\beta_{j+1}}
\vol(R\cap \O^j)^{\alpha_j-\alpha_{j+1}}
\Bigr\}=
\\
=
\prod\limits_{j=1}^n
\frac
{
1-p^{-(\alpha_j-n+j)}}
{(1 -p^{\beta_j+j-1}) (1-p^{-(\alpha_j+\beta_j-n+j)})}
.\end{gather}

{\sc Remark.} The condition of the convergence of this series
is (2.2).

{\bf 2.2. Deformation of $l_2$.}
Fix $\alpha\in \R$.
Let $R,S$ range in $\L_n$. Consider the kernel
\begin{equation}
K_\alpha(R,S)=\frac {\vol(R\cap S)^{\alpha}}
  {\vol(R)^{\alpha/2} \vol(S)^{\alpha/2}}
.\end{equation}

{\sc Theorem 2.2.} a) {\it Let
\begin{equation}
\alpha=0, 1, 2, \dots, n-1 \qquad \text {or $\alpha>n-1$}
.\end{equation}
Then there exist
a  Hilbert space $H_\alpha$
and a system of
 vectors $e_R\in H_\alpha$, where $R$ ranges in $\L_n$,
such that
the scalar products of $e_R$ have the form
\begin{equation}
\langle e_R,e_S\rangle
= K_\alpha(R,S)
\end{equation}
and linear combinations of the vectors $e_R$
are dense in $H_\alpha$
}

\smallskip

b) {\it If $\alpha$ is outside the set {\rm (2.6)},
then a Hilbert space $H_\alpha$
does not exist.}

\smallskip

{\sc Remark.}  The space $H_\alpha$  is unique
in the following sense.
Let $H_\alpha'$ be another
Hilbert space and $e'_R\in H'$ be another
system
of vectors satisfying the same conditions.
Then there exists a unitary operator  $A:H_\alpha\to H_\alpha'$
such that $Ae_R=e_R'$ for all the lattices $R$.

\smallskip

{\sc Remark.}
The space $H_0$ is one-dimensional, all the vectors
$e_R$ in this case coincide.

\smallskip

The group $\GL_n(\K)$ acts in the space $H_\alpha$
by the unitary operators
$$
U_\alpha(g) e_R
=e_{gR}
.$$

{\sc Remark.} Let $\alpha\to+\infty$. Then
$$
\lim\limits_{\alpha\to\infty} K_\alpha(R,T)=
\begin{cases}
1, & \text{if $R=T$;}\\
0, & \text{if $R\ne T$.}
\end{cases}
$$
In this sense, the limit of $H_\alpha$
as $\alpha\to +\infty$ is
the space $l_2(\L_n)$.

\smallskip

{\bf 2.3. Plancherel formula for the spaces $H_\alpha$,
$\alpha>n-1$.}
The following theorem gives the explicit spectral decomposition
of the space $H_\alpha$ with respect to the action
of $\GL_n(\K)$.

\smallskip

{\sc Theorem 2.3.} {\it For $\alpha>n-1$ we have
\begin{equation}
K_\alpha(R,T)=
\int
\Phi_{is_1,\dots, is_n}(R,T)\,d\mu_\alpha(s)
,\end{equation}
where the spherical kernels $\Phi$
were defined in {\rm 1.7},
 the integration is taken over the cube (torus)
$s_j\in [0,2\pi/\ln p]$, and
the Plancherel measure
$\mu_\alpha$ is given by
\begin{multline}
d\mu_\alpha(s)=
C\prod\limits_{j=1}^n (1-p^{-\alpha+j-1})
\prod\limits_{j=1}^n \Bigl| 1-p^{-(\alpha-n+1)/2-is_j}\Bigr|^{-2}
\times \\ \times
\prod\limits_{1\le k<l\le n}
\Bigl|\frac{p^{is_k}-p^{is_l}}
        {p^{is_k}-p^{-1+is_l}}\Bigr|^2
\,ds_1\dots ds_n
,\end{multline}
 where}
$$
C=\frac  {\ln^n p}{ (2\pi)^n n! }
 \prod\limits_{j=1}^n \frac{1-p^{-j}}{1-p^{-1}}
.$$

{\bf 2.4. Plancherel formula for $\alpha=0,1,\dots,n-1$.}

\smallskip

{\sc Theorem 2.4.}  {\it Let $\alpha=0,1,\dots,n-1$. Then
\begin{multline*}
K_\alpha(R,T)=
\int
\Phi_{is_1,\dots, is_\alpha,-(n-\alpha+1)/2, -(n-\alpha+1)/2+1,
  \dots, (n-\alpha+1)/2 }(R,T)\,d\mu_\alpha(s)
,\end{multline*}
 where  $s_j\in[0,2\pi/\ln p]$, the Plancherel measure is
$$
d\mu_\alpha(s)=C_\alpha
\prod\limits_{j=1}^\alpha
   \Bigl|1-p^{(n-\alpha-1)/2+is_j} \Bigr|^2
\prod\limits_{1\le k<l\le \alpha}
\Bigl|\frac{p^{is_k}-p^{is_l}}
        {p^{is_k}-p^{-1+is_l}}\Bigr|^2
\,ds_1\dots ds_\alpha
,$$
and }
$$
C_\alpha=\frac{\ln^\alpha p}{(2\pi)^\alpha (n-\alpha)!}
  (1-p^{-1})^{-\alpha}
\prod\limits_{k=1}^\alpha (1-p^{-k})
\prod\limits_{k=n-\alpha+1}^n (1-p^{-k})
.$$

{\bf 2.5. Spectra.}

\smallskip

{\sc Proposition 2.5.}  {\it Spectrum of each representation
$U_\alpha$ of the group $\GL_n(\K)$ in $H_\alpha$
is multiplicity-free. }

\smallskip

Now we intend to describe the spectra explicitly.

Consider the flag
$$\K^1\subset\K^2\subset\dots\subset \K^k$$
in $\K^n$. Denote by $P_k$ the stabilizer
of this flag in $\GL_n(\K)$. It is the group
of block $(1+1+\dots+1+(n-k))\times(1+1+\dots+1+(n-k))$
upper triangular matrices
$$
A=
\begin{pmatrix}
a_{11}&a_{12}&\dots  & a_{1k} & a_{1(k+1)}     &  \dots  & a_{1n}    \\
 0    &a_{22}&\dots  & a_{2k} & a_{2(k+1)}     &  \dots  & a_{2n}    \\
\vdots&\vdots&\ddots & \vdots & \vdots         &  \ddots & \vdots    \\
0     & 0    &\dots  & a_{kk} & a_{k(k+1)}     &  \dots  & a_{kn}    \\
0     & 0    &\dots  & 0      & a_{(k+1)(k+1)} &  \dots  & a_{(k+1)n}\\
\vdots&\vdots&\ddots & \vdots & \vdots         &  \ddots & \vdots    \\
0     & 0    &\dots  & 0      & a_{n(k+1)}     &  \dots  & a_{nn}    \\
\end{pmatrix}
.$$

In particular, $P_n=B_n$ is the group of all upper triangular matrices.

The spectrum is the support of the Plancherel  measure,
hence
Theorems 2.3--2.4 imply the following corollary.

\smallskip

{\sc Corollary 2.6.}
a) {\it For $\alpha>n-1$ the spectrum consists of
all spherical representations
of $\GL_n(\K)$ of the unitary nondegenerate principal
series}

\smallskip

b) {\it For $\alpha=0,1,2,\dots, n-1$ the spectrum consists
of representations
of $\GL_n(\K)$ unitary induced from the characters
$$\chi_s(A)=\prod_{j=1}^\alpha |a_{jj}|^{is_j},\qquad
0\le s_1\le\dots\le s_\alpha\le 2\pi/\ln p$$
of the group $P_\alpha$.}

\smallskip

{\sc Remarks.}   a) For $\alpha=n-1$ the spectrum
consists of representations $T_{is}$ of unitary nondegenerate
principal series such that $s_n=0$.

b) For $\alpha=n-2,n-3,\dots$ the spectrum consists of
degenerated unitary principal series.
In notation of \S 1.4, consider the subspace $Y_\alpha\subset L^2(\Fl_n)$
of functions that do not depend of the terms
$V_{n-\alpha+1}$,\dots, $V_{n-1}$ of a flag
$\cV\in\Fl_n$; evidently, $Y_\alpha \simeq L^2(\GL_n(\K)/P_k)$.
It is clear that for our values of the parameters
of the representations, the subspace
$Y_\alpha$ is $\GL_n$-invariant. This gives
a realization of our representations in
$L^2(\GL_n(\K)/P_k)$.

\smallskip

{\bf 2.6. A realization of the system $e_R$.}
For definiteness,
let $\alpha>n-1$. Then $H_\alpha$ is equivalent to a direct
integral of representations of principal series.
Here we write explicitly the image of the system $e_R$
in this direct integral.

Let the simplex $\Xi_n$ be the same as above (1.19).
Equip this simplex with the Plancherel measure
$n!\,d\mu_\alpha$, see Theorem 2.3.
Consider the space $\Xi_n\times\Fl_n$
equipped with the product-measure $n!\,d\mu_\alpha\times d\cV$
and  the space $L^2(\Xi_n\times\Fl_n)$
with respect to our measure.
The group $\GL_n(\K)$ acts in the space $L^2$
by the unitary operators $\rho_\alpha(g)$ given by
the formula  (1.20); the formula for the operators themselves
does not contain $\alpha$, but the space of a representation
depends on $\alpha$.

For a lattice $R$,  we consider the function
$u_R(s,\cV)$ on $\Xi_n\times\Fl_n$
given by
$$
u_R(s,\cV)=
\prod\limits_{j=1}^{n-1}
\vol(R\cap V_j)^{1-i s_j+is_{j+1}}
\cdot
\vol(R)^{-(n-1)/2-is_n}
.$$

Our Plancherel formula implies the following corollary.

\smallskip

{\sc Theorem 2.7.}   For each $R$, $T\in\L_n$,
\begin{multline*}
\langle u_R, u_S\rangle_{L^2(\Xi_n\times \Fl_n)}=
\int\limits_{\Xi_n\times\Fl_n}
u_R(s,\cV) \overline {u_T(s,\cV)}\,  n!\,d\mu_\alpha(s)\,d\cV=
\\ =
                      \langle e_R, e_S\rangle_{H_\alpha}
=K_\alpha(R,T)
\end{multline*}

\smallskip

{\sc Remark.} Denote by $\mu_\infty$
Macdonald's Plancherel measure
(1.17)
on $\Xi_n$. The system $u_R(s,\cV)$ is an orthonormal basis
in the space $L^2(\Xi_n\times \Fl_n)$ with respect
to the measure $\mu_\infty\times d\cV$:
$$
\int\limits_{\Xi_n\times\Fl_n}
u_R(s,\cV) \overline {u_T(s,\cV)}\,  n!\,d\mu_\infty(s)\,d\cV=
 \begin{cases} \text{1, if $R=T$}\\
                \text{0, if $R\ne T$}
  \end{cases}
$$

\smallskip

{\bf 2.7. Linear dependence.}   Let
$\alpha=0,1,2,\dots,n-1$.
In these cases,   the space $H_\alpha$ is smaller
than for $\alpha>n-1$ (first, the support of the Plancherel
measure has lesser dimension; second, the representations
in the spectrum are 'smaller').
For these values of $\alpha$ the vectors $e_R$ are linear dependent.

\smallskip

{\sc Theorem 2.8.} {\it
Let $R\subset S$ be arbitrary lattices
such that
$S/R\simeq (\Z/p\Z)^{\alpha+1}$. Then}
\begin{equation}
\sum_{k=0}^{\alpha+1}
\biggl[
(-1)^k
p^{k(k-\alpha-1)/2}
\,\,
\sum\limits_{Q\in \L_n: \,\,
R\subset Q\subset S,\,\, S/Q\simeq (\Z/p\Z)^k}
 e_R\biggr]=0
.\end{equation}

{\bf 2.8. Remark. Relation with the Weil representation.}
Consider the Weil representation (see \cite{Wei1},
see also explicit formulae in \cite{Naz}) of the group
$\Sp(2n,\K)$
and its  restriction $\sigma$  to the subgroup $\GL(n,\K)$.
 The representation
$\sigma$ is equivalent to the natural representation
of $\GL_n(\K)$ in the subspace $L^2_+\subset L^2(\K^n)$
consisting of even functions, i.e., $f(-x)=f(x)$, by the formula
$$
\sigma(g)f(x)= f(gx) |\det g|^{1/2};\qquad g\in\GL(n,\K) ,\quad x\in\K^n
.$$
In particular, the representation $\sigma$ is single-valued
(the Weil representation itself is double-valued).

Denote by  $\xi$ the $\Sp(2n,\O)$-fixed vector in
the Weil representation.
 The 'vacuum vector'
$\xi$ corresponds to the function
$f(z)$ defined by: $f(z)=1$ for $z\in \O^n$ and $f(z)=0$
otherwise.

Consider all the possible vectors $U(g)\xi$, where $g$ ranges
in $\GL_n(\K)$.
We have $U(gh)\xi=U(g)\xi$ for any $h\in\GL_n(\O)$
and hence the vectors of the form $U(g)\xi$
are enumerated by points of $\GL_n(\K)/\GL_n(\O)=\L_n$;
 we denote
the vector corresponding to a lattice $R$ by $\xi_R$.

\smallskip

{\sc Proposition 2.9.}
 $ \langle\xi_R,\xi_S\rangle=K_1(R,S)$.

\smallskip

This is trivial but important. We observe that
$$K_\alpha(R,T)= |\langle\xi_R,\xi_S\rangle|^\alpha$$
In the real case, exactly this procedure gives
the Berezin kernels. Thus, this is a way to conjecture
a form of 'Berezin kernels' in the $p$-adic case.

 The most serious {\tt a posteriori} argument for our analogy
is the formula (2.9) for the Plancherel measure,
since it can be obtained from the formula over $\R$
(see \cite{Ner2}) by replacing
of the usual $\Gamma$-functions by their $p$-adic analogs.

\smallskip

{\bf 2.9. Remark.  Realization of $H_\alpha$ in functions on $\L_n$.}
For each $h\in H_\alpha$, we assign the function $f_h$
on $\L_n$ by the formula
$$
f_h(R)=\langle h,e_R\rangle_{H_\alpha}
.$$
Thus we realize $H_\alpha$ as some space of functions
on $\L_n$. After this, it is possible to apply the usual
formalism of reproducing kernels, see, for instance
\cite{Ner2}, we do not discuss this here.

 For $\alpha=0,1,\dots,n-1$,
 relations (2.10) become difference
equations for functions $f_h$; for each lattices
$R\subset S$ satisfying $S/R=(\Z/p\Z)^{\alpha+1}$,
any function $f_h$ satisfies
$$
\sum_{k=0}^{\alpha+1}
\biggl[
(-1)^k\,p^{k(k-\alpha-1)/2} \,
\sum\limits_{Q\in \L_n: \,\,
R\subset Q\subset S,\,\, S/Q\simeq (\Z/p\Z)^k}
 f_h(R)\biggr]=0
$$

 These equations
 are the analog
of the determinant systems of partial differential equations
defining the degenerated Berezin spaces, see, for instance
\cite{Ner2}.

\smallskip

{\bf 2.10. Remark. Additional symmetry.}
For each element $a=(a_1,\dots,a_n)\in\K^n$
 we define the linear functional
$$\ell_a(z)=a_1z_1+\dots +a_nz_n.$$
Let $R$ be a lattice in $\K^n$.  We define the dual
lattice $R^\vee$ as the set of    all $a\in\K^n$ such that
$\ell_a(z)\in\O$ for all $z\in R$.

\smallskip

{\sc Lemma 2.10.}
$K_\alpha(R^\vee,S^\vee)=K_\alpha(R,S)$

\smallskip

{\sc Proof.}  This follows   from
\begin{align*}
R^\vee\cap S^\vee=(R+S)^\vee,\qquad
\vol(R^\vee)=1/\vol(R),\qquad\qquad\qquad\qquad\\
\qquad \qquad\qquad\qquad
\vol(R+S)\vol (R\cap S)=\vol(R)\cdot\vol(S)
\qquad\qquad\qquad\qquad
\square
\end{align*}

Hence the map $e_R\mapsto e_{R^\vee}$
defines the unitary operator in the space  $H_\alpha$.

Obviously, for $g\in\GL_n(\K)$, $R\in \L_n$, we have
$(gR)^\vee=(g^t)^{-1}R^\vee$,
where ${}^t$ denotes the transposition
of matrices. Thus, we obtain the canonical action
of the semidirect product $(\Z/2\Z)\ltimes\GL_n(\K)$
in the space $H_\alpha$.

\smallskip

{\sc Remark.} In the Weil representation (see above \S2.8),
the operator $e_R\mapsto e_{R^\vee}$ corresponds to the matrix
$\begin{pmatrix}0&1\\-1&0\end{pmatrix}\in{\rm Sp}(2m,\K)$.
In the $L^2(\K^n)$-model of $U_1$
(see \S2.8, or more generally $U_m$, see below \S6.1) it corresponds
to the Fourier transform.

  Theorem 2.8 and  the structure
of the Plancherel  formula in Theorem 4.4
(see below)
also are consistent with this symmetry.

\smallskip

{\bf 2.11. Remark. Identities with Hall--Littlewood functions.}
Let $p$ be real, $p>1$. Define the functions
$\phi_\lambda^{(p)}(\kk)=\phi_\lambda(\kk)$ by the formula (1.11).
For prime values of $p$ these functions are spherical functions
over $\K=\Q_p$.

\smallskip

{\sc Proposition 2.11.} a) {\it The identity {\rm (2.8)} holds
for arbitrary $p>1$, $\alpha>n-1$.}

\smallskip

b) {\it For arbitrary $p>1$ and
 sufficiently large real $\alpha$,
$$\sum\limits_{\kk}
\phi_{is}^{(p)}(\kk) p^{-\alpha\sum_j |k_j|} \nu(\kk)
=\prod\limits_{l=0}^{n-1} (1-p^{-\alpha+l})
\prod\limits_{j=1}^n \bigl|1-p^{-(\alpha+n-1)/2+is_j}\bigr|^{-2}
,$$
where $\nu(\kk)$ is given by} (1.4).

\smallskip

{\bf 2.12. Remark: the Tamagawa zeta-function.}
Assume in Theorem 2.1
$$\alpha_j=t,\qquad \beta_j=\gamma_j-t.$$
Passing to the limit as $t\to+\infty$,
we obtain

{\sc Corollary 2.12.}
$$
\sum\limits_{R\in\L_n,R\subset\O^n}
\prod\limits_{j=1}^n \vol(R\cap \K_j)^{\gamma_j-\gamma_{j+1}}=
\prod\limits_{j=1}^n (1-p^{-\gamma_j-n-j})^{-1}
.$$
This is the Tamagawa zeta-function, see \cite{Mac2}, IV.4.

\smallskip

{\bf 2.13. Remark. Lattices in $\Q^n$ and classical zeta-function.}
Consider the $n$-dimensional space $\Q^n$ over the
rational numbers. A {\it lattice} in $\Q^n$ is a free
$\Z$-submodule with $n$ generators. We denote the space of lattices
in $\Q^n$ by $\L_n(\Q)$.

For a lattice $R\subset\Q^n$, we denote by
$\upsilon_n(R)$ the volume of the torus $\R^n/R$.

Denote by $e_1$, \dots,  $e_n$ the standard basis
in $\Q^n$. Denote
$$
\Q^k=\Q e_1\oplus\dots\oplus \Q e_k,
\qquad \Z^k=\Z e_1\oplus\dots\oplus\Z e_k
.$$
Fix complex $\alpha_1$, \dots, $\alpha_n$
and $\beta_1$, \dots, $\beta_n$; assume
$\alpha_{n+1}=0$, $\beta_{n+1}=0$.
Let
$$
\Re\alpha_j-n+1>1,\qquad
\Re\alpha_j+\Re\beta_j-n+j>1,\qquad
\Re\beta_j+j<0
.
$$

{\sc Theorem 2.13.}
\begin{multline}
\sum\limits_{R\in\L_n(\Q)}\,\,
\prod\limits_{j=1}^n
\upsilon_k(R\cap\Q^k)^{-\beta_k+\beta_{k+1}}
\upsilon_k(R\cap\Z^k)^{-\alpha_k+\alpha_{k+1}}
=\\=
\prod\limits_{j=1}^n
\frac{\zeta(-(\beta_j+j-1))\zeta(\alpha_j+\beta_j-n+j)}
{\zeta(\alpha_j-n+j)}
,\end{multline}
{\it where
$$
\zeta(s)=\sum\limits_{m=1}^\infty
  \frac 1{m^s}=
\prod\limits_{\text{$p$ is prime}}
\Bigl(1-\frac 1{p^s}\Bigr)^{-1}
$$
is the zeta-function.}

\bigskip

{\bf \large 3. Calculation of the beta-sum}

\nopagebreak

\medskip

\addtocounter{sec}{1}
\setcounter{equation}{0}

Here we prove Theorems 2.1 and 2.13.

\smallskip

{\bf 3.1. Transformation of the beta-sum.}
Denote by $\L_{n-1}$ the space of lattices in the subspace
$\K^{n-1}\subset\K^n$.
We transform the expression (2.3) to the form
\begin{gather}
\sum_{R'\in \L_{n-1}}
\prod\limits_{j=1}^{n-2}\Biggr[
\Bigl\{
\vol(R'\cap \K^j)^{\beta_j-\beta_{j+1}}
\vol(R'\cap \O^j)^{\alpha_j-\alpha_{j+1}}
\Bigr\}\times \\ \times
\vol(R'\cap \K^{n-1})^{\beta_{n-1}}
\vol(R'\cap \O^{n-1})^{\alpha_{n-1}}\times\\
\times \sum\limits_{R\in\L_n:\, R \cap \K^{n-1}=R'}
\Bigl\{
\frac{\vol(R)^{\beta_n}}{\vol(R\cap \K^{n-1})^{\beta_n}}
\frac{\vol(R\cap \O^n)^{\alpha_n}}{\vol(R\cap \O^{n-1})^{\alpha_n}}
\Bigr\}\Biggr]
.\end{gather}

{\sc Lemma 3.1.}
{\it The interior sum {\rm (3.3)} is equal to
$$\frac {\sigma_n}{\vol(\O^{n-1}\cap R')},$$
where}
\begin{equation}
\sigma_n= \frac {1-p^{-\alpha_n}}
    {(1-p^{-\alpha_n-\beta_n})(1-p^{\beta_n+n-1})}
.\end{equation}

Lemma 3.1 will be proved below.

\smallskip

{\sc Corollary 3.2.} {\it
Denote by $F_n(\alpha;\beta)$ the left  hand side of
{\rm (2.3)--(2.4)}. The following identity
for $F_n(\alpha;\beta)$ holds
$$
F_n(\alpha_1,\dots,\alpha_n; \beta_1,\dots, \beta_n)=
\sigma_n
F_{n-1}(\alpha_1,\dots,\alpha_{n-1}-1; \beta_1,\dots, \beta_{n-1})
, $$
 where $\sigma_n$ is given by} (3.4).

This corollary implies Theorem 2.1.

\smallskip

{\bf 3.2. New variables $\xi$ and $\nu$.}
Consider the subspace $\K^{n-1}\subset \K^n$.
We have $\K^n=\K^{n-1}\oplus \K e_n$.

Fix a lattice $R'\subset \K^{n-1}$.
Let $R$ satisfy the property
\begin{equation}
R\cap \K^{n-1}=R'.
\end{equation}

Each vector $r$ in $R$ has the form
$r=a\, p^k e_n+v$, where $v\in\K^{n-1}$,
$|a|=1$.

\smallskip

{\sc Lemma 3.3.} {\it Denote by $\xi$ the minimal possible $k$,
let $h=p^\xi e_n+v\in R$. Then }
\begin{equation}
R=R'\oplus \O h=R'\oplus \O(p^\xi e_n +v).
\end{equation}

\smallskip

This is obvious. The next Lemma is equivalent to
Lemma 3.3.

\smallskip

{\sc Lemma 3.$3'$.} {\it
Consider the natural map $\pi:\K^n\to \K e_n$.
Let  $h\in R$ and $\pi(h)$ generates
the $\O$-module  $\pi(R)$.
Then $R=R'\oplus \O h$.}


 \smallskip

{\sc Lemma 3.4.} {\it Lattices $R_1=R'\oplus \O (p^{\xi_1} e_n+v_1)$
and $R_2=R'\oplus\O (p^{\xi_2} e_n+v_2)$ coincide iff
$\xi_1=\xi_2$, $v_1-v_2\in R'$.}

\smallskip

{\sc Proof} is obvious.

\smallskip

Thus  a lattice $R$  with the property
(3.5) is determined by the integer $\xi=\xi(R)$ and
 an element $v=v(R)$ of the quotient group $\K^{n-1}/R'$.

\smallskip

We also define a nonnegative integer $\nu=\nu(R)$
by the condition

\smallskip

 --- if $v\in R'+\O^{n-1}$, then $\nu=0$;

\smallskip

 --- otherwise, $\nu$ is the minimal $k$ such that
                      $p^k v\in R'+\O^{n-1}$.

\smallskip

{\sc Lemma 3.5.}
\begin{gather*}
a)\qquad \frac{\vol(R)}{\vol(R')}=p^{-\xi}.\\
b) \qquad \frac{\vol(R\cap \O^n)}{\vol(R\cap \O^{n-1})}=
p^{-\max(0,\xi+\nu)}
.\end{gather*}

{\sc Proof.}
The statement a) is obvious.
The statement b) is a corollary of the following lemma.

\smallskip

{\sc Lemma 3.6.}
{\it Let $R$ has the form {\rm (3.6)}. Then  $R\cap\O^n$
has the form
$$
R\cap\O^n=R'\oplus \O\bigl[p^\kappa(p^\xi e_n+v)+w\bigr]
,$$
where $\kappa=\max (\nu,-\xi)$ and $w\in-p^\kappa v + R'$.}

\smallskip

{\sc Proof.} First,
\begin{equation}
p^\kappa(p^\xi e_n+v)+w  \in  R'\qquad \Longleftrightarrow          \qquad
  \kappa\ge 0,\quad w\in R'
\end{equation}
Secondly,
\begin{equation}
p^\kappa(p^\xi e_n+v)+w  \in  \O^n\qquad\Longleftrightarrow  \qquad
  \kappa+\xi\ge 0,\quad p^\kappa v+w\in \O^{n-1}
\end{equation}
Since $w\in R'$, $p^\kappa v+w\in\O^{n-1}$,
 we have $p^\kappa v\in R'+\O^{n-1}$.
But $\kappa\ge 0$, hence $\kappa\ge \nu$.
Also, $\kappa\ge-\xi$, see (3.8).
It remains to apply Lemma 3.3-$3.3'$.
\hfill $\square$

\smallskip

%
%
%
%
\smallskip

{\sc Lemma 3.7.}
{\it The number of lattices $R\in \L_n$ satisfying (3.5)
with given  $\eta\ge 0$, $\xi$  is }
$$
u_n(\nu,\xi)=\begin{cases}
\frac {\textstyle1}{\textstyle\vol (R'\cap \O^{n-1})},& \nu=0; \\
   \textstyle \frac
{\textstyle p^{\nu(n-1)} (1-p^{-n+1})}
{\textstyle\vol (R'\cap \O^{n-1})},& \nu>0
\end{cases}
$$

{\sc Proof.} Nothing depends on $\xi$.
We must find number of vectors
$v\in\K^{n-1}/R'$ such that $p^\nu v\in (R'+\O^{n-1})/R'$.

The index of $R'$ in $R'+\O^{n-1}$ is
$1/\vol(R'\cap \O^{n-1})$.
It remains to find the
number $u^*_n(\nu)$ of $v\in\K^{n-1}/(R'+\O^{n-1})$ such that
$p^\nu \in R'+\O^{n-1}$, since we have
$$u_n(\nu,\xi)=
\frac {u^*_n(\nu)}{\vol(R'\cap \O^{n-1})}
.$$

For each lattice $S\subset\K^{n-1}$, the additive group
$\K^{n-1}/S$
is isomorphic to $ (\K/\O)^{n-1}$.
Hence we must find the number of solutions
of the  equation  $p^\nu v=0$ in $ (\K/\O)^{n-1}$.

The additive group $\K/\O$ is the inductive limit of the cyclic
groups
$$
\K/\O=\lim_{k\to\infty}\Z/ p^k \Z
.$$
Hence the number of solutions of the equation $p^\nu v=0$
in  $ (\K/\O)^{n-1}$ is $p^{\nu(n-1)}$.
But $p^{(\nu-1)(n-1)}$ of these solutions are also
solutions of the equation $p^{\nu-1}v=0$.

If $\nu=0$, we have a unique solution $v=0$.
\hfill $\square$

\smallskip

{\bf 3.3. The end of calculation.}
Lemma 3.5, 3.7 reduce  an evaluation of (3.3)
to an evaluation of the sum
$$
\frac 1{\vol(R'\cap \O^{n-1})}
\sum_{\xi\in \Z}\sum_{\nu\ge 0}
p^{-\alpha_n \max(0,\xi+\nu)} p^{-\beta_n\xi} u_n(\nu,\xi)
$$
Here the summation amounts to
summations of several geometric progressions.

\smallskip

{\bf 3.4. Proof of Theorem 2.13.}
Let $\Q_p$ be $p$-adic numbers, $\O_p$ be $p$-adic integers.
For each prime $p$ consider the natural embeddings
$\Q\to\Q_p$, $\Q^n\to\Q^n_p$ and the induced map
$$
\pi_p:\L_n(\Q)\to \L_n(\Q_p)
,$$
i.e., for $R\in\L_n(\Q)$ we consider its closure in $\Q_p^n$.

Our calculation is based on the following two remarks.

First, for any $R\in\L_n(\Q)$,
$$
\upsilon_n(R)^{-1}=\prod\limits_{\text{$p$ is prime}}
\vol_{\Q^n_p}\bigl( \pi_p(R)\bigr)
.$$

Secondly,  the map
$$
R \mapsto
\bigl(\pi_2(R), \pi_3(R), \pi_5(R), \pi_7(R), \pi_{11}(R),\dots\bigr)
$$
is a bijection of $\L_n(\Q)$ and the set of sequences
$(S_2,S_3,S_5,S_7,S_{11},\dots)$ such that
$S_p\in\L_n(\Q_p)$ and $S_p=\O^n_p$ except finite
number  of $p$; see \cite{Wei}, Theorem V.2.2,

\smallskip

Hence the left hand side of (2.11) transforms to the form
$$
\prod\limits_{\text{$p$ is prime}}
\Bigl[
\sum_{S_p\in \L_n(\Q_p)}
\prod\limits_{j=1}^n
\Bigl\{
\vol(S_p\cap \Q^j_p)^{\beta_j-\beta_{j+1}}
\vol(S_p\cap \O^j_p)^{\alpha_j-\alpha_{j+1}}
\Bigr\}
\Bigr]
.$$

 It remains to apply Theorem 2.1.

\medskip

{\bf\large 4. Plancherel formula}

\nopagebreak

\medskip

\addtocounter{sec}{1}
\setcounter{equation}{0}

This section contains proofs of Theorem 2.3
(in \S 4.1), Theorem 2.4 (in \S\S 4.2--4.6),
Proposition 2.5 (in \S 4.7), Theorem 2.7
(in \S 4.8), Proposition 2.11 (in \S 4.9).
We also obtain 'indefinite Plancherel formula'
(Theorem 4.4), which is used below in \S5.

{\bf 4.1. Plancherel formula.}
Consider the $\GL_n(\O)$-invariant  function
$$
\Delta_\alpha(R)=K_\alpha(\O^n,R)=\frac{\vol(R\cap \O^n)^\alpha}
           {\vol (R)^{\alpha/2}} ,
\qquad R\in\L_n
.$$

By (1.15),
its spherical transform is
\begin{equation}
\widetilde\Delta_\alpha(\lambda)=
\sum_{R\in\L_n}
\frac{\vol(R\cap \O^n)^\alpha}
  {\vol (R)^{\alpha/2}}
   \cdot
    \prod\limits_{j=1}^{n-1}
\cdot
\vol\bigl(R\cap \K^j\bigr)^{1+\lambda_j-\lambda_{j+1}}
\cdot\vol(R)^{-(n-1)/2+\lambda_n}
.\end{equation}

By  Theorem 2.1,
it equals
$$
\widetilde\Delta_\alpha(\lambda)=
\prod_{j=1}^n \frac
{1-p^{-(\alpha-n+j)})}
{(1-p^{-(\alpha-n+1)/2+\lambda_j})
 (1-p^{-(\alpha-n+1)/2-\lambda_j)})  }
.$$
Applying the inversion formula (1.16) for
the spherical transform, we
obtain Theorem 2.3.

\smallskip

{\bf 4.2. Plancherel formula as a contour integral.}
Introduce the new variables $z_k=p^{is_k}$.
Also  introduce a new notation for the spherical functions
$$
\phi[z_1,\dots,z_n;R]:=\phi_{is_1,\dots,is_n}(R)
.$$

Now Theorem 2.3 is converted to the form
\begin{equation}
\Delta_\alpha(R)=C_0 I_0,
\end{equation}
where  $I_0$ is the following integral over a torus
\begin{align}
\!\!\!\!\!\!\!\!\!\!\!\!\!\!\!\!\!\!\!\!\!\!
I_0=I_0(\alpha;R)=
\!\!\!\!\!\!\!
\int\limits_{|z_1|=1,\dots,|z_n|=1}
\prod_{j=1}^n
\frac 1
 {(z_j-p^{(\alpha-n+1)/2}) (z_j-p^{-(\alpha-n+1)/2})}
\times
\\
\times
\prod_{1\le k<l\le n} \frac{(z_k-z_l)^2}
      {  (z_k-pz_l)   (z_k-p^{-1}z_l) }
\,
\phi[z_1,\dots,z_n;R ]\, dz_1\dots dz_n
\end{align}
and the constant $C_0$ is
\begin{equation}
C_0=C_0(\alpha)=\frac{p^{\alpha n}}{(-2\pi i)^n}\,\,
\prod_{j=1}^n
\frac{1-p^{-j}}{1-p^{-1}}
\,\,
\prod_{l=0}^{n-1}
(1-p^{-\alpha+l})
\end{equation}

The  function $\Delta_\alpha(R)$ is a holomorphic function
in  $\alpha\in\C$.
The integral expression $C_0 I_0$ for  $\Delta_\alpha(R)$
is holomorphic in the domain $\Re \alpha>n-1$.
We intend to obtain the holomorphic continuation
of $C_0(\alpha)I_0(\alpha;R)$ into the
whole plane $\alpha\in\C$.

\smallskip

{\bf 4.3. Analytic continuation  to the strip
$n-3<\Re \alpha<n-1$.}
Denote
\begin{equation}
\beta:=(\alpha-n+1)/2,
\end{equation}
\begin{equation}
\mu_m(z_1,\dots,z_m):=
\prod_{1\le k<l\le m} \frac{(z_k-z_l)^2}
      {  (z_k-pz_l)   (z_k-p^{-1}z_l) }
\end{equation}

It can be easily checked that $\mu_m$ is symmetric with respect to
$z_j$. Recall also that the spherical functions
are symmetric with respect to $z_j$.

\smallskip

\begin{figure}
\1
\end{figure}

{\sc Lemma 4.1.} {\it For
$n-3<\Re \alpha<n-1$,
\begin{equation}
\Delta_\alpha(R)=C_0\, I_0\, +\, n\, C_{01}\,( I^+_{01}\, +
                      \, I^-_{01})   \,
                   + \, n(n-1)\, C_1\, I_1
,
\end{equation}
where $C_0$, $I_0$ are the same as above {\rm (4.3)--(4.5)},
\begin{align}
I^+_{01}=
\int\limits_{|z_1|=1,\dots,|z_{n-1}|=1}
\prod_{j=1}^{n-1}
\frac{z_j-p^\beta}
{(z_j-p^{\beta+1})(z_j-p^{\beta-1})(z_j-p^{-\beta})}
\times \\ \times
\mu_{n-1}(z_1,\dots,z_{n-1})\,\phi[z_1,\dots,z_{n-1}, p^\beta;R]\,
dz_1\dots dz_{n-1}
\notag
,\end{align}

\begin{align}
I^-_{01}=
\int\limits_{|z_1|=1,\dots,|z_{n-1}|=1}
\prod_{j=1}^{n-1}
\frac{z_j-p^{-\beta}}
{(z_j-p^{-\beta-1})(z_j-p^{-\beta+1})(z_j-p^{\beta})}
\times \\ \times
\mu_{n-1}(z_1,\dots,z_{n-1})\,\phi[z_1,\dots,z_{n-1}, p^{-\beta};R]\,
dz_1\dots dz_{n-1}
\notag
,\end{align}

\begin{align}
I_1=
\int\limits_{|z_1|=1,\dots,|z_{n-2}|=1}
\prod_{j=1}^{n-2}
\frac{(z_j-p^\beta)(z_j-p^{-\beta})}
{(z_j-p^{-\beta-1}) (z_j-p^{\beta+1})
(z_j-p^{\beta-1})  (z_j-p^{-\beta+1}) }
\times\\ \times
\mu_{n-2}(z_1,\dots,z_{n-2})\,
\phi[z_1,\dots,z_{n-2}, p^{-\beta},  p^{\beta};R]\,
dz_1\dots dz_{n-2}
\notag
,\end{align}
and the constants are given by}
$$
C_{01}=\frac{2\pi i C_0}{p^\beta-p^{-\beta}}
,$$
$$
C_1=p^{-2\beta} (1-p^{-2\beta-1})^{-1}(1-p^{-2\beta+1})^{-1} (2\pi i)^2 C_0
.$$
We denote the integrands in $I_k$,  $I_{kl}^\pm$,
see (4.9)--(4.11) and Theorem 4.4 below
by
$$
\fI_k,\,\,\,   \fI_{kl}^\pm
$$

{\sc Proof.} First, we expand
\begin{equation}
\prod_{j=1}^n
\frac 1
 {(z_j-p^{\beta}) (z_j-p^{-\beta})}=
(p^\beta-p^{-\beta})^{-n}
\prod_{j=1}^n
\Bigl( \frac 1{z_j-p^\beta} - \frac 1{z_j-p^{-\beta}} \Bigr)
\end{equation}
and open brackets. The integral (4.3)--(4.4) splits
into the sum of $2^n$ integrals and it is sufficient
to construct  the analytic continuation of each summand.

For definiteness,  consider the summand
\begin{multline*}
N(\beta):=\int\limits_{|z_1|=1,\dots,|z_{n}|=1}
\frac 1
{(z_1-p^\beta)\dots(z_s-p^\beta)
(z_{s+1}-p^{-\beta})\dots  (z_n-p^{-\beta})}
\times\\ \times
\mu_{n}(z_1,\dots,z_{n})\,\phi[z_1,\dots,z_{n};R]\,
dz_1\dots dz_{n}
;\end{multline*}
by symmetry considerations, we do not lose a generality.
Denote the integrand in $N(\beta)$ by $\fN(\beta;z)$.
The analytic  continuation  of $N(\beta)$
through the line $\Re \beta=0$
is given by the same integral  over
another contour
\begin{equation}
|z_1|=1-\epsilon,\, \dots,\, |z_s|=1-\epsilon,\,
|z_{s+1}|=1+\epsilon,\, \dots,\, |z_n|= 1+\epsilon
.\end{equation}
Consider the family of contours $L_0$, \dots $L_n$,
where $L_m$ is given by
$$
L_m:\qquad
|z_j|=
\begin{cases}
1, & \text{if $j\le m$,}\\
1\pm \epsilon  &\text{if $j>m$,}
\end{cases}
,$$
and the signs $\pm$ are the same as in (4.13).
In particular,
$L_n$ is the torus $|z_1|=\dots=|z_n|=1$,
and $L_0$ is (4.13).

Then
\begin{equation}
\int_{L_0} =\Bigl(\int_{L_0}- \int_{L_1}\Bigr) +
 \Bigl(\int_{L_1}- \int_{L_2}\Bigr) +\dots
\Bigl(\int_{L_{n-1}}- \int_{L_n}\Bigr)+ \int_{L_n}
.\end{equation}
Each bracket can be evaluated by one-dimensional
residues.
In the last bracket we obtain
$$
\int\nolimits_{|z_1|=\dots=|z_{n-1}|=1}
\res_{z_n=p^{-\beta}}
\fN(\beta;z)\prod\limits_{j=1}^{n-1} dz_j
,
$$
here we obtain a desired expression.

For other brackets the  picture is more complicated.
For instance, in the first bracket we obtain
$$
-\int\res_{z_1=p^{\beta}}
\fN(\beta;z)\prod\limits_{j=2}^{n} dz_j
,$$
 where the integration is taken
over the torus $|z_2|=1-\epsilon$, \dots,
 $|z_n|=1+\epsilon$. For this torus we apply the
transformation of (4.14) type,
etc., etc. Thus we obtain $2^n$ summands having the
form $\int \res\res\dots\res$. However, $\mu_n(z)=0$
on the hyperplanes $z_k=z_l$, and hence
$$\res\limits_{z_k=p^\beta} \res\limits_{z_l=p^\beta}\fN(\beta,z)=0,
\qquad
\res\limits_{z_k=p^{-\beta}} \res\limits_{z_l=p^{-\beta}}\fN(\beta,z)=0
$$
Hence
$$
\int_{L_0}-\int_{L_n}=
-\sum_{j\le s} \int \res\limits_{z_j=p^\beta}\fN
 +
\sum_{j> s} \int \res\limits_{z_j=p^{-\beta}}\fN
-
\sum_{j\le s, k>s}
\int \res\limits_{z_j=p^\beta}\res\limits_{z_k=p^{-\beta}}
\fN
,$$
 where integrations are taken over tori $|z_l|=1$.
Next, we transform this sum to the form
\begin{multline}
-\sum_{1\le j\le n}\int \res\limits_{z_j=p^\beta}
\fN(\beta;z)\!\!\!\prod\limits_{1\le m\le n; m\ne j}\!\!\! dz_m
 +
\sum_{1\le k\le n} \int \res\limits_{z_j=p^{-\beta}}
\fN(\beta;z)\!\!\!\prod\limits_{1\le m\le n; m\ne k}\!\!\! dz_m
                                   -  \\
\sum_{1\le j\le n, 1\le k\le n}
\int \res\limits_{z_j=p^\beta}\res\limits_{z_k=p^{-\beta}}
\fN(\beta;z)\prod\limits_{1\le m\le n; m\ne j,k} dz_m
\end{multline}
(all the 'new' summands equal zero).

Adding together  all the $2^n$ summands of the integral
$I_0$, we obtain an expression of the form (4.15), only
$\fN$ is replaced by the integrand $\fI_0$ of $I_0$.
Evaluatiing the residues and
applying the symmetry
of the spherical functions with respect to $z_j$,
we obtain (4.8);
the symmetry with respect to $z_j$
 is the origin of the coefficients $n$,  $n(n-1)$ in
(4.8).

\smallskip

{\bf 4.4. Plancherel formula at $\alpha=n-1$ and $\alpha=n-2$.}

\nopagebreak

\smallskip

{\sc Proposition 4.2.}
$$
\Delta_{n-2}(R)=
n(n-1)\,\bigl[ C_1\, I_1 \bigr]\Bigr|_{\alpha=n-2}
.$$

{\sc Proof.} For $\alpha=n-2$ we have
 $C_0=C_{01}=0$ in (4.8).
This gives the required result.

\smallskip

\begin{figure}
\2
\end{figure}

{\sc Proposition 4.3.}
$$
\Delta_{n-1}(R)=\frac 12 n\,
\bigl[ C_{01}\,( I^+_{01}\, +
                       \, I^-_{01})\bigr] \Bigr|_{\alpha=n-1}
.$$

{\sc Proof.}
For $\alpha=n-1$,  the both formulae (4.2) (4.8), are not
valid. Consider the new function
$$h(\beta;R):=
\frac{2\beta\ln p}{C_0(\beta)}
 \Delta_\alpha(R)
$$
(where $(\alpha-n+1)/2=\beta$). The factor
$2\beta\ln p/C_0(\beta)$ is a holomorphic   nonvanishing
function near the point $\beta=0$  and hence
it is sufficient to find an integral expansion for $h(0,R)$.
We have
$$
h(0,R)=\frac{1}{2\pi i}
\int\limits_{|\beta|=\rho} \frac{h(\beta,R)}{\beta}\,d\beta
$$
for each sufficiently small $\rho$.
We transform this expression to the form
\begin{multline}
h(0,R)=
\frac{2\ln p}{2\pi i}
\int\limits_{|\beta|=\rho}
I_0(\beta)\,d\beta
+\\+
\frac{2\ln p}{2\pi i}
\int\limits_{|\beta|=\rho,\,\,\Re\beta<0}
\Bigl[
\frac{n \cdot 2\pi i}{p^\beta-p^{-\beta}}
\Bigl(I^+_{01}(\beta) + I^-_{01}(\beta)\Bigr)
+  \\ +
\frac{n(n-1)(2\pi i)^2}{2p^{2\beta}(1-p^{-2\beta-1})(1-p^{-2\beta+1})}
I_1(\beta)\Bigr] \,d\beta
.\end{multline}

First, the function $I_0(\beta)$ is even, i.e.,
$I_0(\beta) =I_0(-\beta)$, see (4.3).
Hence, the first term in (4.16) vanishes.

Secondly, $I_1(0)$ is finite,
and the corresponding scalar factor  is finite
at $\beta=0$. Hence the last term of (4.16) tends to
zero as $\rho\to 0$.

Thirdly, let us evaluate
\begin{equation}
\lim\limits_{\rho\to 0}
\int\limits_{|\beta|=\rho, \,\Re \beta<0}
\frac{2\ln p}{p^\beta-p^{-\beta}}
I_{01}^+(\beta;R)\,d\beta
\end{equation}
In formula (4.9), $I_{01}^+$
 was defined as an integral over the contour
$|z_1=\dots=|z_{n-1}=1$. We can replace this contour
by $M:\, |z_1|=\dots=|z_{n-1}|=1+\epsilon $, and
$\int_M \fI_{01}^+$ is holomorphic in the strip
$|\Re\beta| <\epsilon$.
This allows to transform (4.17)
to the form
$$
\int\limits_M
\lim\limits_{\rho\to 0}    \!\!\!\!\!\!
\int\limits_{|\beta|=\rho, \,\Re \beta<0} \!\!\!
\frac{2\ln p}{p^\beta-p^{-\beta}}
\,\fI_{01}^+ (\beta,z)\,d\beta\,\, dz_1\dots dz_{n-1} =
\pi i
\int\limits_M
\,\fI_{01}^+ (0,z)\,d\beta\,\, dz_1\dots dz_{n-1}
.$$
This implies the required result.


\smallskip

{\bf 4.5. Complete analytic continuation.}
We preserve the notation (4.6)--(4.7)
for $\beta$ and $\mu(z)$

\smallskip

{\sc Theorem 4.4.}
{\it Let $\Re\alpha\ne n-1, n-3,\dots, -n+1$.
 Then
\begin{multline}
\Delta_\alpha(R)= C_0 I_0+
\sum\limits_{k:\, k\ge 1,\, 2(k-1)+\Re\alpha<n-1}
\frac{n!}{(n-2k)!} C_k I_k+\\+
\sum\limits_{k,l:\, k\ge 0,\, l\ge 1,\,2(k+l-1)+\Re\alpha<n-1 }
\frac{n!}{(n-2k-l)!}\,C_{kl}\, (I_{kl}^++I_{kl}^-)
,
\end{multline}
where
\begin{multline*}
I_k=          \!\!\!\!\!\!
\int\limits_{|z_1|=\dots=|z_{n-2k}|=1}
\prod\limits_{j=1}^{n-2k}
\frac{(z_j-p^{-\beta-k+1})   (z_j-p^{\beta+k-1})}
{(z_j-p^{\beta+k})   (z_j-p^{-\beta-k})
 (z_j-p^{\beta-1})   (z_j-p^{-\beta+1})}
\times \\ \times
\phi[z_1,\dots,z_{n-2k},
p^{\beta} , p^{\beta+1},\dots, p^{\beta+k-1}
 p^{-\beta-k+1},\dots, p^{-\beta+1},p^{-\beta};R]
\times \\ \times
\mu_{n-2k}(z_1,\dots,z_{n-2k}) \,dz_1\dots dz_{n-2k}
,\end{multline*}
\begin{multline*}
I_{kl}^+=          \!\!\!\!\!\!
\int\limits_{|z_1|=\dots=|z_{n-2k-l}|=1}
\prod\limits_{j=1}^{n-2k-l}
\frac{(z_j-p^{-\beta-k+1})   (z_j-p^{\beta+k+l-1})}
{(z_j-p^{\beta+k+l})   (z_j-p^{-\beta-k})
 (z_j-p^{\beta-1})   (z_j-p^{-\beta+1})}
\times \\ \times
\phi[z_1,\dots,z_{n-2k-l},
p^{\beta} , p^{\beta+1},\dots, p^{\beta+k+l-1}
p^{-\beta-k-1},\dots, p^{-\beta+1}, p^{-\beta};R]
\times \\ \times
\mu_{n-2k-l}(z_1,\dots,z_{n-2k-l}) \,dz_1\dots dz_{n-2k-l}
\end{multline*}
\begin{multline*}
I_{kl}^-=          \!\!\!\!\!\!
\int\limits_{|z_1|=\dots=|z_{n-2k-l}|=1}
\prod\limits_{j=1}^{n-2k-l}
\frac{(z_j-p^{\beta+k-1})   (z_j-p^{-\beta-k-l+1})}
{(z_j-p^{-\beta-k-l})   (z_j-p^{\beta+k})
 (z_j-p^{\beta-1})   (z_j-p^{-\beta+1})}
\times \\ \times
\phi[z_1,\dots,z_{n-2k-l},
p^{\beta} , p^{\beta+1},\dots, p^{\beta+k-1}
p^{-\beta-k-l-1},\dots, p^{-\beta+1}, p^{-\beta};R]
\times \\ \times
\mu_{n-2k-l}(z_1,\dots,z_{n-2k-l}) \,dz_1\dots dz_{n-2k-l}
,\end{multline*}
and the constants $C_k$, $C_{kl}$ are given by}
\begin{equation}
C_k=(2\pi i)^{2k} C_0\cdot\frac{p^{2\beta k+k(k-1)} (1-p)^{2k} (1-p^{2\beta-1})}
                 {1-p^{2\beta+2k-1}}
\cdot\prod\limits_{j=1}^k
  \frac 1
  {(1-p^j)^2(1-p^{2\beta+j-2})^2}
,\end{equation}
\begin{multline}
C_{kl}= (2\pi i)^l
C_k\cdot\frac{p^{-(\beta+k)l-l(l-1)/2}(1-p^{-2\beta-2k+1})(1-p^{-1})^l}
             {1-p^{-2\beta-2k-l+1}}
\times\\ \times
\prod_{j=1}^l\frac 1
    {(1-p^{-k-j}) (1-p^{-2\beta-k-j+2})}
.\end{multline}

{\sc Proof.} We start from Lemma 4.1
and intend to write the analytic continuation
of (4.8) to the strip $|\Re \beta+1|<\epsilon$.

First, $I_0$ has no singularities   in the half-plane
$\Re\beta<0$.
Thus we must write the analytic continuation of $I_{01}^\pm$, $I_1$,
see (4.9)--(4.11), to the strip $|\Re \beta+1|<\epsilon$.

\smallskip

Secondly, the integrand $\fI_1$ of $I_1$ contains the factor
$$
\prod \frac 1
  {(z-p^{\beta+1})    (z-p^{-\beta-1}) }
.$$
We transform this factor as (4.12), and repeat literally
the proof of Lemma 4.1.
This gives 3 new summands $I_2$, $I_{11}^+$, $I_{11}^-$
in the strip
$-\epsilon<\Re\beta+1<0$.

Thirdly,  $\fI_{01}^+$ contains the factor
$$
\prod \frac 1
  {(z-p^{\beta+1})}
;$$
 all other factors of $\fI_{01}^+$
 have no poles for $\beta$ lying in our
strip  $|\Re \beta+1|<\epsilon$. Hence we write the analytic
continuation of $I_{01}^+$  as
$$
\int\nolimits_{|z_1|=1-\epsilon,\dots, |z_{n-1}|=1-\epsilon}
\fI_{01}^+\, dz_1\dots dz_n
,$$
 and evaluate the analytic
continuation as the sum of residues
$$
\int\limits_{|z_1|=\dots=|z_{n-1}|=1} \fI_{01}^+ \prod\limits_{j=1}^{n-1} dz_j-
\sum\limits_{j=1}^{n-1}
\int\limits_{|z_1|=\dots=|z_{j-1}|=|z_{j+1}|=|z_{n-1}| }
\res\limits_{z_j=p^{\beta+1}} \fI_{01}^+
\prod\limits_{1\le j\le n-1,\, j\ne s} dz_j
$$
(by the symmetry considerations all the summands of the sum
$\sum_{j=1}^{n-1}$ coincide).
Other summands
$$
\int \res\limits_{z_j=p^{\beta+1}} \res\limits_{z_m=p^{\beta+1}} \fI_{01}^+,
\qquad
\int \res\limits_{z_j=p^{\beta+1}} \res\limits_{z_m=p^{\beta+1}}
     \res\limits_{z_s=p^{\beta+1}}\fI_{01}^+,\qquad \text{etc.}
$$
vanish, since the factor $\mu_{n-1}(z)$ is zero
at the hyperplanes $z_k=z_m$.
This gives the new summand $I_{02}^+$
in the strip
$-\epsilon<\Re\beta-1<0$.

The case of $I_{01}^-$ is similar.

\smallskip

Further, $I_1$, $I_{01}^\pm$ have no singularities
in the half-plane
$\Re\beta<-1$.
The 'new' summands $I_2$, $I_{11}^\pm$, $I_{02}^\pm$
have singularities on the line $\Re\beta=-2$.
Their analytic continuation through this line
can be obtained  in the same way, etc., etc.

Now, our explicit  formulae for $I_k$, $I_{kl}^\pm$,
$C_k$, $C_{kl}^\pm$
 can be proved by induction.
In fact it is necessary to check the identities
$$
C_{k+1}\fI_{k+1}=2\pi iC_{k}  \!\!\!
\res\limits_{z_{n-2k}=p^\beta}
\res\limits_{z_{n-2k-1}=p^{-\beta}}
\fI_{k};
\quad
C_{k(l+1)}\fI_{k(l+1)}^\pm= 2\pi i
C_{kl}\!\!\!
\res\limits_{z_{n-2k-l}=p^{\mp\beta}}
\fI_{kl}^\pm
$$

\smallskip

{\bf 4.6. Integer values of $\alpha$.
 Proof of Theorem 2.4.}
Let now $\alpha=0,1,2,\dots,n-1$.  In this case,
the factor (see (4.5))
\begin{equation}
\prod_{m=0}^{n-1} (1-p^{-\alpha+m})
\end{equation}
of $C_0$ is decisive.
It  vanishes at all our $\alpha$, and hence the summand
$I_0$ in the Plancherel formula (4.18) disappears.
The same factor kills the
most of other summands, since
$C_0$ is present as a factor in $C_k$ and $C_{kl}$.

\smallskip

{\sc Theorem 4.5.} {\it Let $n-\alpha=2k$ be positive
even integer. Then}
$$
\Delta_{n-2k}= \frac{n!}{(n-2k)!}
\bigl[C_k I_k\bigr]\Bigr|_{\alpha=n-2k}
$$

{\sc Proof.}
All other summands of the formula (4.18) vanish due (4.21), the summand
$C_k I_k$ survives, since the denominator of $C_k$
contains the factor
$$1-p^{2\beta+2k-1}=1-p^{\alpha-n+2k}$$
(see (4.19)), which also vanishes at $\alpha$.

\smallskip

{\sc Theorem 4.6.} {\it Let $n-\alpha=2m+1$ be positive
odd integer. Then}
$$
\Delta_{n-2m-1}=\frac{n!}{2\, (n-2m-1)!}
\bigl[ C_{m1}(I_{m1}^+(\alpha) + I_{m1}^-(\alpha))\bigr]
\Bigr|_{\alpha=n-2m-1}
$$

{\sc Proof.} Theorem 4.4 does not give
immediate answer in this case and we use
the same arguments as in Proposition 4.3.

Consider the analytic continuation of our integral
to the strip $-m<\Re \beta<-m+1$, or equivalently
$n-2m-1<\Re \alpha< n-2m+1$, it is given by (4.18).
We must evaluate the limit of each summand of (4.18)
as $\beta$ tends to $-m$ from our strip.

First, let $k<m$. Then $I_k$ is holomorphic in
$\Re \beta<-m+1$, and $C_k\bigr|_{\beta=-m}=0$. Hence, the
 summand $C_k I_k$  vanishes at $\beta=-m$.

Secondly, let $k+l<m$. Then $I_{kl}^\pm$
are holomorphic in $\Re \beta<-m+1$, and
$C_{kl}\bigr|_{\beta=-m}=0$. Hence $C_{kl} I_{kl}^\pm$
vanish.

Thirdly, let $k+l=m$, $l>0$. Then the analytic continuation
of $I_{kl}^+$ to the strip $|\Re\beta+m|<\epsilon$ is given by
$$
\int\nolimits_{|z_1|=\dots=|z_{n-2k-l}|=1\mp\epsilon}
\fI^+_{kl}
$$
 Hence $I_{kl}^+$
is holomorphic in  $|\Re\beta+m|<\epsilon$. But the coefficient
$C_{kl}$ vanishes at $\beta=-m$ again.

Thus the problem is reduced to the evaluation
of $\lim_{\beta\to -m+0} C_m I_m$.
The factor $C_m$ has a simple zero at $\beta=-m$,
and hence it is sufficient to evaluate
$$\lim_{\beta\to -m+0} (\beta+m) I_m.$$

We represent
\begin{multline}
\frac{(z-p^{\beta+m-1}) (z-p^{-\beta-m+1})}
 {(z-p^{\beta+m}) (z-p^{-\beta-m})
  (z-p^{\beta-1}) (z-p^{-\beta+1})}
 =\\=
\frac{\lambda(\beta)}
 { (z-p^{\beta+m}) (z-p^{-\beta-m})}
+
\frac{\sigma(\beta)}
{(z-p^{\beta-1}) (z-p^{-\beta+1}) }
,\end{multline}
where
$$
\lambda(\beta)=
\frac{(1-p)(1-p^{2\beta+2m-1})}
     {(1-p^{m+1})(1-p^{2\beta+m-1})};\qquad
\sigma(\beta)=\frac{p(p^{2\beta+k-2}-1)(p^k-1)}
     {(p^{2\beta+k-1} -1)(p^{k+1}-1)}
.$$

Let
$J$ be a subset in $\{1,2,\dots,n-2m\}$. Denote by $\|J\|$
the number of elements of  $J$.
According (4.22), we represent $(\beta+m)I_m$ as a sum of $2^{n-2m}$
integrals
\begin{multline}
(\beta+m)I_m(\beta)=
\sum\limits_{J}
\lambda(\beta)^{\|J\|}
\sigma(\beta)^{n-2m-\|J\|}
\int\limits_{|z_s|=1,\, s\notin J}
\prod\limits_{s\notin J} \frac 1{(z_s-p^{\beta-1})  (z_s-p^{-\beta+1})}
\times\\ \times
\Biggl[    (\beta+m)
\int\limits_{|z_j|=1,\, j\in J}
\prod\limits_{j\in J} \frac 1{(z_j-p^{\beta+m})  (z_j-p^{-\beta-m})}
\times\\ \times
\mu_{n-2m}(z_1,\dots,z_{n-2m})\phi(\dots;R)
\prod\limits_{j\in J} dz_j \Biggr]
\prod\limits_{s\notin J} dz_s
.\end{multline}

Denote by $F(\beta)=F_J(\beta)$ the expression in big brackets
(it depends also on  $z_s$ for $s\notin J$).
Denote by $\fF$ the integrand in $F(\beta)$.
The analytic continuation of $F(\beta)$
can be easily written as in
proof of Lemma 4.1. In the strip $-m-1<\Re \beta<-m$ it is
$$F(\beta)+(\beta+m)G(\beta),$$
 where $F(\beta)$ is the same expression
in square brackets (and hence $F(\beta)$ is singular
on the line $\Re\beta=-m $) and
\begin{multline*}
G(\beta) =
2\pi i
\sum\limits_{k\in J}\int\limits_{\text{$|z_j|=1$ for $j\in J$, $j\ne k$}}
\Bigl[-\res\limits_{z_k=p^{\beta+m}} \fF+
     \res\limits_{z_k=p^{-\beta-m}} \fF \Bigr]
\prod\limits_{j\in J,\, j\ne k}\,dz_j
-\\-
 (2\pi i)^2
\sum\limits_{k,l\in J}
 \int\limits_{\text{$|z_j|=1$ for $j\in J$,  $j\ne k,l$}}
\res\limits_{z_k=p^{\beta+m}}
\res\limits_{z_l=p^{-\beta-m}}\fF
\prod\limits_{j\in J,\, j\ne k,l}\,dz_j
\end{multline*}

Thus,
$$  \lim\limits_{\beta\to -m+0}
F(\beta)=
\frac 1{2\pi i} \int\nolimits_{|\beta+m|=\rho}
  \frac {F(\beta)\,d\beta}{\beta+m}
+
\frac 1{2\pi i} \int\nolimits_{|\beta+m|=\rho, \Re\beta<-m}
 G(\beta)\,d\beta
$$
The transformation
$\beta\mapsto -2m-\beta$
preserves  $F_I(\beta)$
and hence the first summand in right hand side is 0.

The second summand gives the required result
as in proof of Proposition  4.3.

\smallskip

{\bf 4.7. Absence of multiplicity.}
Evidently, the $\GL_n(\O)$-fixed vector $e_{\O^n}$
is $\GL_n(\K)$-cyclic. This easily implies
Proposition 2.5, see proof of Lemma 1.10 in \cite{Ner2}.


\smallskip

{\bf 4.8. System $e_R$ in the direct integral.
Proof of Theorem 2.7.}
First, $\rho_\alpha(g)u_R=u_{gR}$.
Hence  we can assume $R=\O^n$, and
$u_{\O^n}(s,\cV)=1$.
Thus, we must evaluate
$$\int_{\Xi_n}\int_{\Fl_n} u_T(s,\cV)\,d\cV \,d\mu_\alpha(s)$$
 Integrating over $\Fl_n$ using (1.8),
we obtain
$$\int_{\Xi_n} \phi_{is}(T) \, \,d\mu_\alpha(s)$$
It equals $K_\alpha(\O^n,T)$ by Theorem 2.3.

\smallskip

{\bf 4.9. Identities with Hall--Littlewood functions.}
Let as prove the statement a) of Proposition 2.11.
Using (1.11), it is possible to  represent the integral (4.3)--(4.4)
as a sum of residues. The result is a (nonhand)
 rational
expression in $p$ and $p^\alpha$.
Thus the required identity  is an identity of rational
functions which holds for countable number of $p$
(and for all $\alpha$ if $p$ is fixed).
Hence the identity holds always.

The statement b) follows from a) and the inversion formula
for spherical transform
(the Macdonald's proof of Theorem 5.1.2
is based on calculations with Hall--Littlewood-type
expression (1.11) and do not use primality of $p$).

\medskip

{\bf\large 5. Positivity and nonpositivity}

\medskip

\addtocounter{sec}{1}
\setcounter{equation}{0}

Here we prove Theorem 2.2 on positive definiteness
of the kernel $K_\alpha$.

\smallskip

{\bf 5.1. The case
 $\alpha>n-1$.} This follows from
Theorem 2.3, since $K_\alpha(R,T)$
is an integral of positive definite kernels
$\Phi_{is_1,\dots,is_n}(R,T)$
with a positive weight $d\mu_\alpha$.

\smallskip

{\bf 5.2. Positive definiteness of the kernel $K_\alpha(R,T)$
for integer $\alpha$.}
This again follows from the Plancherel formula (Theorem 2.4).

More simple  proof is given
below in 6.1.

\smallskip

{\bf 5.3. The case $\alpha<0$.}
This case is obvious:
$$\langle e_R,e_R\rangle=\langle e_S,e_S\rangle=1,\qquad
\langle e_R,e_S\rangle>1
.$$
This is impossible.

\smallskip

{\bf 5.4. Noninteger $\alpha$ between 0 and $n-1$.}
Assume that the function $\Delta_\alpha$ is positive definite.
Then it can be expanded as an integral of positive definite
spherical functions  $\phi_\tau$
\begin{equation}
\Delta_\alpha(R)=\int \phi_\tau(R) d\kappa_\alpha(\tau)
\end{equation}
with respect to some positive measure $\kappa_\alpha$.

In Theorem 4.4, we  obtained some expansion
\begin{equation}
\Delta_\alpha(R)=\int \phi_\tau(R) d\sigma_\alpha(\tau)
\end{equation}
 of $\Delta_\alpha$
as an integral of spherical functions
with respect to
 some measure (charge) $\sigma_\alpha$.
For all noninteger $\alpha<n-1$, the summand $I_{01}$
is present in our expansion with a nonzero coefficient.
But the spherical function $\phi_{is_1,\dots,is_{n-1}, n-1-\alpha}$
is not positive definite.
Hence the expansions (5.1), (5.2) are different.

 It is sufficient
to reduce the existence of these two different expansions
to a contradiction.   We do this in the rest of this section.

\smallskip

{\bf 5.5. The set $\Sigma_n$.} Denote by $\rho$
the vector $(n-1)/2,(n-1)/2-1, \dots, -(n-1)/2\in\R^n$.
The symmetric group $S_n$ acts on $\R^n$ by
permutations of the coordinates. Denote by $Q$ the convex hull
of $S_n$-orbit of the vector $\rho$.

By $\Sigma_n$ we denote the set of all
$$\tau=(\tau_1,\dots,\tau_n)\in\C^n/\tfrac{2\pi i}{\ln p}\,\Z^n$$
such that each $\tau_j$ is real or pure imaginary
and $\Re \tau\in Q$.

The set $\Sigma_n$ is important for us by the following reasons.

\smallskip

$1^\circ.$ If $\phi_\tau$ is positive definite, then $\tau\in \Sigma_n$.
In particular, the measure $\kappa_\alpha$
from (5.1) is supported by
$\Sigma_n$

\smallskip

$2^\circ.$ For $\alpha>0$, the measure $\sigma_\alpha$
from (5.2) is supported by  $\Sigma_n$.

$3^\circ.$ $\Sigma_n$ is compact.

\smallskip

The unitary dual of $\GL_n(\K)$ is known (see \cite{Tad1}, \cite{Tad2}),
but it is easier to reduce the statement $1^\circ$ to
 Theorem 4.7.1 of \cite{Mac1}.

We also denote by $\Sigma^u_n$ the set of positive definite
spherical functions.

\smallskip

{\bf 5.6. Action of the Hecke algebra.}
The Hecke algebra $\cH_n$
  of $\GL_n(\K)$ (see \cite{Mac2})
is the convolution algebra of compactly supported
 $\GL_n(\O)$-{\it biinvariant} functions on
$\GL_n(\K)$.

A function $f$ on $\GL_n(\K)$ is called $\GL_n(\O)$-biinvariant
if
$$f(h_1gh_2)=f(g)\qquad \text{for $g\in\GL_n(\K)$, $h_1,h_2\in\GL_n(\O)$}
.$$
The multiplication in the Hecke algebra is the convolution.

The Hecke algebra acts on the space of $\GL_n(\O)$-biinvariant
functions by the convolutions.

The spherical functions $\phi_\tau$ are the eigenfunctions
of $\cH_n$, i.e.,
\begin{equation}
\gamma*\phi_\tau=c_\gamma(\tau) \phi_\tau, \qquad \gamma\in\cH
,
\end{equation}
where  $c_\gamma(\tau)$ is a scalar factor.

Denote by $\cA$ the algebra of all polynomial
expressions of
$p^{\tau_1}$, \dots, $p^{\tau_n}$, $p^{-\tau_1}$, \dots, $p^{-\tau_n}$
symmetric with respect to permutations of $\tau_j$.
The map $\gamma\mapsto c_\gamma$ is an isomorphism
$\cH\to\cA$, see \cite{Mac2}, V.3.2.

\smallskip

{\sc Lemma 5.1.} {\it Let a positive definite
$\GL_n(\O)$-biinvariant function
$f$   be represented as an integral over $\Sigma_n$ with respect
to some charge
$$
f(g)=\int_{\Sigma_n} \phi_\tau(g)\, d\sigma(\tau)
.$$
Let $u(\tau)\in\cA$ be   nonnegative
on $\Sigma_n$. Then  the function
$$
q(g)=\int_{\Sigma_n} u(\tau) \phi_\tau(g)\, d\sigma(\tau)
.$$
also is positive definite.}

\smallskip

{\sc Proof.} Since $f$ is positive definite, it admits an expansion
$$
f(g)=\int_{\Sigma_n^u}
\phi_{\tau}(g) d\nu(\tau)
.$$
for some positive measure $\nu$.
Let $\gamma$ be the element of the Hecke
algebra corresponding to the function $u(\tau)$.
Then
$$(\gamma*f)(g)=\int_{\Sigma_n} u(\tau) \phi_\tau(g)\, d\sigma(\tau)
=\int_{\Sigma_n^u} u(\tau) \phi_\tau(g)\, d\nu(\tau)$$

The second integral defines a positive definite function,
and this proves Lemma.  \hfill $\square$

\smallskip

{\bf 5.7. Approximation of spherical functions by matrix elements
of $U_\alpha$.}
Denote by $W_n$ the hyperoctahedral group of $\R^n$, i.e.,
the group generated by permutations
of the coordinates and arbitrary changes of directions
of axes.

Consider   the map    $\Sigma_n\to \R^n$   given by
$$
(\tau_1,\dots,\tau_n)\mapsto
  (p^{\tau_1}+p^{-\tau_1},\dots,p^{\tau_n}+p^{-\tau_n})
$$
  Obviously, it is constant on
 orbits of $W_n$ and moreover this map is an embedding of the quotient
space $\Sigma_n/W_n$ to $\R^n$.

Consider the functions
$$
u_k(\tau)=\sum_{j=1}^n   \bigl(p^{\tau_j}+ p^{-\tau_j}\bigr)^k
.$$
Obviously, these functions  are real on the set $\Sigma_n$
and they separate orbits of the hyperoctahedral group
on $\Sigma_n$.

\smallskip

Fix a point $\kappa=(is_1,\dots,is_{n-1},n-1-\alpha)\in\Sigma_n$.
Consider the function
$$F(\tau)=\sum_{k=1}^n \bigl[u_k(\tau)-u_k(\kappa)\bigr]^2. $$
The function $F$ is zero at the orbit $W_n\cdot\kappa$
and is positive outside this orbit.
Let $M$ be the maximum of $F$ on $\sigma$.
Consider the new function
$$X(\tau)=M-F(\tau).$$
It is a positive function having a maximum
on the orbit $W_n\cdot \kappa$.
For simplicity, assume that $s_j$ are pairwise distinct.
Then the points $w\tau$ are nondegenerate critical points
of $X$.
Let $\xi$ be the corresponding element
 of the Hecke algebra.
Consider the sequence   of positive definite functions
\begin{equation}
m^{(n-1)/2}\xi^m*\Delta_\alpha(g)=
m^{(n-1)/2} \int X^m(\tau) \phi_s(g) d\sigma_\alpha(\tau)
.\end{equation}
It can be easily checked, that its limit
as $m\to\infty$
is
\begin{equation}
r(g):=\sum_{w\in W_n} c(w\kappa) \phi_{w\kappa}(g)
,\end{equation}
 where the scalar coefficients $c(w\kappa)$
are nonzero.
This function is positive definite as a limit
of positive definite functions. Hence $r(g)$
can be  expanded into an integral of positive definite
spherical functions,
\begin{equation}
r(g)=\int_{\Sigma_n^u} \phi_\tau(g) \,d\lambda(\tau)
.\end{equation}
Let us evaluate the limit of $\xi^m*r$ as $m\to\infty$
 (now without the normalizing factor
$m^{(n-1)/2}$).
Applying it to (5.6), we obtain
\begin{equation}
\xi^m*r(g)=\int_{\Sigma_n^u} X(g)^m \phi_\tau(g) \,d\lambda(\tau)
.\end{equation}
The orbit $W_n\chi$ has no intersection with the (closed) set
$\Sigma_n^u$, and hence the limit of integrals (5.7)
is zero.
If we evaluate the same limit using (5.5), we  obtain $r(g)$
(since $\xi^m*r=r$).
Hence $r(g)=0$. But $r(g)$ is a finite linear combination
of the spherical functions, on the other hand the spherical functions
are the eigenfunctions of the Hecke algebra.
This is a contradiction.


\medskip

{\bf\large 6. Integer values of $\alpha$ and difference equations}

\medskip

\addtocounter{sec}{1}
\setcounter{equation}{0}

Here we prove Theorem 2.8 on linear dependence
of the vectors $e_R$.

\smallskip

{\bf 6.1. Embedding of $H_m$ to $L^2(\K^{nm})$.}
Let $m$ be a nonnegative integer.
Consider the space
$$
\K^{mn}=\underbrace{\K^n\oplus\dots\oplus\K^n}_{\text{$m$ times}}
.$$
For a lattice $\R\in\L_n$, we define the lattice
$$
R^m:=R\oplus \dots\oplus R \subset \K^{nm}
.$$
Define the function $e_R$ on $\K^{nm}$ by
the rule
$$
e_R(\theta)=
\begin{cases} \vol (R)^{-m/2},
      & \text{if $\theta\in R^m$}\\
 0, &\text{if $\theta\notin R^m$}
\end{cases}
.$$

Obviously, the $L^2(\K^{nm})$-scalar products
of the vectors $e_R$ are
$$\langle e_R, e_S \rangle=K_\alpha(R,S)$$
and hence we can identify $H_m$ with
the subspace in $L^2(\K^{mn})$ generated by the functions
$e_R(\theta)$.

In particular, this proves the existence of
the Hilbert space $H_m$ for integer $m$.

\smallskip

{\bf 6.2. Reduction of Theorem 2.8 to combinatorial
problem.} Now, let $m=0,1,\dots,n-1$.

Consider lattices $R\subset S$  in $\K^n$ such that
$S/R=(\Z/p\Z)^{m+1}$.   Without loss of generality,
we can assume $\vol(S)=1$.

 It is convenient to identify
the quotient $S/R$ with $(m+1)$-dimensional linear space
$\F_p^{m+1}$ over the $p$-element field $\F_p$.
Denote by $\pi$ the natural projection
$$\pi: \,S\to S/R\simeq \F_p^{m+1}$$

Let $Q\in\L_n$ satisfy $R\subset Q\subset S$.

The function $e_Q$ is zero outside $S^n$;
if $\xi-\eta\in R^n$, then $e_Q(\xi)=e_Q(\eta)$.
This allows to consider  the functions $ e_Q$
as functions
on the quotient group
$$
S^m/R^m\simeq
\underbrace{\F_p^{m+1}\oplus\dots\oplus \F_p^{m+1}}_%
{\text{$m$ times}}\simeq \F_p^{(m+1)m}
.$$

More precisely, for
each
 linear subspace $L\subset\F_p^{m+1}$
we define the function
$$
\wt e_L(w_1\oplus\dots\oplus w_{m}):=
\begin{cases}
p^{\codim(L) m/2} , &\text{ if $w_j\in L$ for all $j$;}\\
0, &\text{otherwise}
.\end{cases}
$$

If $Q$ lies between $R$ and $S$, then $L:=Q/R$ is a
linear subspace in $S/R=\F_p^{m+1}$. Obviously, we have
$$
e_Q(\theta_1\oplus\dots\oplus\theta_m)=
\wt e_L(\pi (\theta_1)\oplus\dots\oplus\pi(\theta_m))
$$
Thus, it is sufficient to investigate the linear dependencies
of the functions $\wt e_L$.

\smallskip

{\bf 5.3. Proof of Theorem 2.8.}
For $k=0,\dots,m+1$ consider the function $G_k$ on
$\F_p^{(m+1)m}$   given by
$$
G_k(w_1\oplus\dots\oplus w_{m})=
p^{-k m/2}
  \sum\limits_{k:\codim L=k} \wt e_L(w_1\oplus\dots\oplus w_{m})
.$$
We intend to find numbers $u_0$, \dots, $u_m$ such that

$$
\sum_k u_k G_k=0
.$$

{\sc Lemma 6.1} $u_k=(-1)^k p^{k(k-1)/2}$.

\smallskip

{\sc Proof.}
Obviously, $G_k(w_1\oplus\dots\oplus w_{m})$
coincides with the number of subspaces $L$ of
codimension $k$ containing all the vectors
$w_1$,\dots,$w_{m}$. Denote by $W$
the linear span of $w_j$.  We must count
subspaces in $\F^{m+1}_p$ containing $W$, or
equivalently, linear subspaces in $\F^{m+1}_p/W$.

The number $A_l^j$ of $j$-dimensional subspaces
in $\F_p^l$ is
$$
A_l^j=
\frac
{(p^l-1)(p^{l-1}-1)\dots(p^{l-j+1}-1)}
{(p^j-1)(p^{j-1}-1)\dots(p-1)}
$$

Hence, our Lemma is equivalent to the family of the identities
\begin{equation}
\sum\limits_{i=0}^{s}
(-1)^i p^{i(i-1)/2}
\frac{(p^s-1)(p^{s-1})\dots(p^{s-i+1}-1)}
{(p^i-1)(p^{i-1}-1)\dots (p-1)}=0
\end{equation}

For this identity. we can refer to \cite{Mac2}, Ex.I.2.3, or to
 the $q$-binomial theorem (see, for instance \cite{GR}, \S 1.3)
\begin{multline*}
\sum\limits_{i=0}^\infty
\frac
{(1-aq)(1-aq^2)\dots(1-aq^{i-1})}
{(1-q)(1-q^2)\dots(1-q^i)} z^i=
\frac{(1-az)(1-azq)(1-azq^2)\dots}
{(1-z)(1-zq^2)(1-zq^3)\dots}
,\end{multline*}
we substitute $q=1/p$, $a=p^s$, $z=q$.

\tt

Math. Phys. Group,

Institute for Theoretical and Experimental Physics,

B. Cheremushkinskaya, 25,  Moscow -- 117259, Russia

neretin@mccme.ru, \qquad neretin@gate.itep.ru

\end{document}